\documentclass[10pt,draft]{amsart}
\usepackage{amssymb}
\usepackage{latexsym}
\usepackage{amsthm}
\usepackage{amsfonts}
\usepackage{amscd}
\usepackage{amsxtra}

\newcommand{\N}{\mathbb{N}}

\providecommand{\gather}{\begin{gather}}
\providecommand{\endgather}{\end{gather}}
\providecommand{\align}{\begin{align}}
\providecommand{\endalign}{\end{align}}
\providecommand{\aligned}{\begin{aligned}}
\providecommand{\endaligned}{\end{aligned}}

\newtheorem{theorem}{Theorem}[section]
\newtheorem{lemma}[theorem]{Lemma}

\newtheorem{definition}[theorem]{Definition}
\newtheorem{proposition}[theorem]{Proposition}
\newtheorem{remark}[theorem]{Remark}
\newtheorem{corollary}[theorem]{Corollary}
\newtheorem{example}[theorem]{Example}

\title[Disjoint and Simultaneously Hypercyclic Pseudo-Shifts]{Disjoint and Simultaneously Hypercyclic Pseudo-Shifts}

\author[N. \c{C}olako\u{g}lu]{Nurhan \c{C}olako\u{g}lu}
\address[N. \c{C}olako\u{g}lu]{Department of Mathematics, Istanbul Technical University, Maslak, 34469, Istanbul, Turkey}
\email{colakn@itu.edu.tr}

\author[\"O. Martin]{\"Ozg\"ur Martin}
\address[\"O. Martin]{Department of Mathematics, Mimar Sinan Fine Arts University, Silah\c{s}\"or Cad. 71, Bomonti \c{S}i\c{s}li 34380, Istanbul, Turkey}
\email{ozgur.martin@msgsu.edu.tr}

\author[R. Sanders]{Rebecca Sanders}
\address[R. Sanders]{Department of Mathematics, Statistics, and Computer Science, Marquette University, Milwaukee, WI 53201}
\email{rebecca.sanders@marquette.edu}

\date{December 09, 2021 }

\begin{document}

\maketitle

\begin{abstract}
	We  characterize disjoint and simultaneously hypercyclic tuples of unilateral pseudo-shift operators on $\ell^p(\N)$. As a consequence, complementing the results of Bernal and Jung, we give a characterization for simultaneously hypercyclic tuples of unilateral weighted shifts. We also give characterizations for unilateral pseudo-shifts that satisfy the Disjoint and Simultaneous Hypercyclicity Criterions. Contrary to the disjoint hypercyclicity case, tuples of weighted shifts turn out to be simultaneously hypercyclic if and only if they satisfy the Simultaneous Hypercyclicity Criterion.
\end{abstract}

\section[Introduction]{Introduction}

Weighted shift operators have always been a rich source of examples and counterexamples in operator theory and, in particular, in linear dynamics which is an area that study the dynamical properties of continuous linear operators on topological vector spaces. Unfortunately, when it comes to joint dynamical notions as disjointness and simultaneousness for finitely many hypercyclic operators, the dynamics of weighted shifts turns out to be very limited. For example, it was shown in \cite{BeMaSa} that tuples of weighted shifts can never be disjoint weakly mixing or disjoint mixing, nor satisfy the Disjoint Hypercyclicity Criterion. The aim of this paper is to study the joint dynamics of pseudo-shifts which is a generalization of weighted shifts in order to provide a new source of examples, unify the known results related to the shift operators, and provide new results for simultaneous hypercyclicity.

Let $\mathbb{N}$ denote the set of positive integers, $X$ be a separable and infinite dimensional Banach space over the  real or complex scalar field $\mathbb{K}$, and let $B(X)$ denote the algebra of bounded linear operators on $X$. Recall that, a bounded linear operator $T$ acting on a Banach space $X$ is called {\it hypercyclic} if there exists a vector $x$ in $X$ such that its orbit $\{T^nx: n \geq 0\}$ is dense in $X$. Such $x$ is called as a {\it hypercyclic vector} for $T$. 

First example of a hypercyclic operator on a Banach space was given by Rolewicz \cite{R} by showing that if $B$ is the unweighted unilateral backward shift on $\ell^(\N)$ with $1 \le p < \infty$, then $\lambda B$ is hypercyclic if and only if $|\lambda| > 1$.  Recall that the unweighted unilateral backward shift $B$ is defined as $Be_m = e_{m-1}$, for $m \geq 2$ and $Be_1 = 0$ where $\{e_m: m \in \mathbb{N} \}$ is the canonical basis of the sequence space $\ell^p(\mathbb{N})$. 

A natural generalization of the unweighted unilateral backward shift $B$ is the unilateral weighted shifts which multiplies the shifted vector with a weight sequence $\omega = (w_n)_{n \in \N}$ of scalars, that is,  $B_{\omega} e_m = w_m e_{m-1}$, for $m \geq 2$ and $B_{\omega} e_1 = 0$. In a fundamental paper in the area,  Salas  \cite{Sa} completely characterized the hypercyclic unilateral weighted backward shifts in terms of their weight sequences on $\ell^p(\mathbb{N})$  with $1 \le p<\infty$.

Recently, joint dynamics of finitely many operators have attracted attention of researchers. In 2007, the notion of disjointness in hypercyclicity is introduced independently by Bernal \cite{BG2} and by B\`es and Peris \cite{BeP1}.

\begin{definition}{\rm
Operators $T_1, \ldots,T_N  \in B(X)$ with $N \ge 2$ are called \textit{disjoint hypercyclic}, or \textit{d-hypercyclic} for short, if the direct sum operator $\oplus_{i=1}^{N} T_i = T_1 \oplus \dots \oplus T_N$ has a hypercyclic vector of the form $(x,\dots, x) \in \oplus_{i=1}^{N}X$.  Such a vector $x$ is called a
\textit{d-hypercyclic vector} for the operators $T_1, \ldots,T_N$.  If the set of d-hypercyclic vectors of $T_1, \ldots,T_N$ is dense in $X$, then the operators $T_1, \ldots,T_N$ are called  \textit{densely d-hypercyclic}. 
}
\end{definition}

The similar and weaker notion of simultaneous hypercyclicity is introduced and studied by Bernal and Jung \cite{BGJ} in 2018.

\begin{definition} \cite[Definition 2.1]{BGJ}\label{D:s-hyper} {\rm
	Operators $T_1, \ldots, T_N \in B(X)$ with $N \ge 2$ are called \textit{simultaneously hypercyclic}, or \textit{s-hypercyclic} for short, if there exists a vector $x \in X$ such that 
	\[
		\overline{\{(T_{1}^n x,\ldots,T_{N}^n x):n \in \mathbb{N} \} }
		\supset\Delta\left( \oplus_{i=1}^{N} X \right),
	\]
	where $\Delta\left( \oplus_{i=1}^{N} X \right) = \left\{(x,\ldots,x): x \in X \right\}$ denotes the diagonal of $\oplus_{i=1}^{N} X$. Such a vector $x$ is said to be a \textit{s-hypercyclic vector} of $T_1, \ldots, T_N$. If the set of s-hypercyclic vectors of the operators $T_1, \ldots, T_N$ is dense in $X$, then the operators $T_1, \ldots, T_N$ are called as \textit{densely s-hypercyclic}.
}
\end{definition}

In both definitions, we assume $\oplus_{i=1}^{N} X$ is endowed with the product topology.

B\`es and Peris \cite{BeP1} and Bernal and Jung \cite{BGJ} have characterized d-hypercyclicity and s-hypercyclicity of tuples of distinct powers of weighted shifts, respectively. They showed that given weighted shifts $B_{\omega^{(1)}}, \ldots, B_{\omega^{(N)}}$ and distict integers $1 \leq r_1 < \ldots < r_N$, d-hypercyclicity and s-hypercyclicty of $B_{\omega^{(1)}}^{r_1}, \ldots, B_{\omega^{(N)}}^{r_N}$ are equivalent. They also showed that these operators are d-hypercyclic (or s-hypercyclic) if and only if they satisfy the d-Hypercyclicity Criterion if and only if they satisfy the s-Hypercyclicity Criterion (see Definition \ref{def:CoHC} and \ref{def:sHC}, respectively). 

Characterization of d-hypercyclicity and s-hypercyclicty of weighted shifts that are raised to the same power turned out to be more delicate. To the best of our knowledge, the problem of characterization of s-hypercyclicity of a finite tuple of weighted shifts (raised to the same power) is still open. We will answer this question in Section 3. In \cite{BeMaSa}, B\`es and the last two authors characterized d-hypercyclicity of weighted shifts. In fact, this work highlighted many differences between the dynamics of one operator and the disjoint dynamics of a finite tuple. In \cite{BeMaSa}, it is shown that weighted shifts can never be d-weakly mixing or d-mixing, nor satisfy the d-Hypercyclicity Criterion (see  \cite{BeMaSa} for the definitions and for more on d-hypercyclic weighted shifts), thus disjoint dynamics of weighted shifts are limited. This is the reason why we turn our attention to a larger family of operators which are still easy to handle but exhibit a larger spectrum of joint chaotic dynamical behavior. 

In \cite{KGGro}, Grosse-Erdmann defined a generalization of weighted shift operators which are called as \textit{pseudo-shifts}, and studied their dynamics. In \cite{BDP}, Bongiorno, Darji, and Di Piazza  studied the dynamics of pseudo-shifts with constant weight sequences and strictly increasing shifting functions which they call as the \textit{Rolewicz-type operators}. This paper led the first two authors to characterize disjoint and simultaneously hypercyclic Rolewicz-type operators in \cite{CoMa}. 

Recently, Wang and Zhou \cite{WaZh18}, and Wang and Liang \cite{WaLi19}, respectively, characterized d-hypercyclicity and d-supercyclicity of tuples of pseudo-shifts of the form $T_{f,\omega^{(1)}}, \ldots, T_{f,\omega^{(N)}}$ which have the same shifting function. In \cite{WaChZh}, Wang, Chen and Zhou characterized d-hypercyclicity and d-supercyclicity of tuples of pseudo-shifts of the form $T^{r_1}_{f_1,\omega^{(1)}}, \ldots, T^{r_N}_{f_N,\omega^{(N)}}$ where powers are pairwise distinct. To the best of our knowledge, a complete study of d-hypercyclic or s-hypercyclic pseudo-shifts raised to the same power hasn't been done until this paper.

Following Bongiorno et al. \cite{BDP}, in Definition \ref{def:pseudoshift}, we slightly change the definition of pseudo-shifts, assuming that the shifting function is strictly increasing but without any restrictions on the weights. This change provides simplicity in statements and proofs.

In Section 2, we give a sufficient condition for satisfying a disjoint blow up/collapse condition and use it to characterize disjoint hypercyclic unilateral pseudo-shifts in terms of their weight sequences and shifting functions, unifying the results in \cite{BeMaSa} and  \cite{BeP1}. In Section 3, we first provide a  Simultaneous Blow-Up/Collapse Criterion, corresponding to the one in Section 2, and use it to characterize simultaneous hypercyclic pseudo-shifts. As an application, we give a characterization for simultaneous hypercyclic unilateral weighted shifts. We show that, contrary to the disjointness case, unilateral weighted shifts are s-hypercyclic if and only if they satisfy the s-Hypercyclicity Criterion. In contrast to this result, we give an example of a tuple of pseudo-shifts which satisfy the Simultaneous Blow-Up/Collapse Criterion but fail to satisfy the s-Hypercyclicity Criterion.

In the rest of the Introduction, we will provide the needed definitions and propositions from \cite{BeP1} and \cite{BGJ}.  
\begin{definition}\label{def:ToplogicalTransitivity}
    \normalfont{We say the operators $T_1, \dots, T_N \in B(X)$ with $N \ge 2$ are \textit{d-topologically transitive} (respectively, \textit{s-topologically transitive}) provided that for any non-empty open sets $V_0, U_1, \dots, U_N$ (respectively, non-empty open sets $U,V$) of $X$, there exists an integer $n \in \mathbb{N}$ such that 
    $V_0 \cap T_1^{-n}(U_1) \cap \cdots \cap T_N^{-n}(U_N) \ne \emptyset$ (respectively, $V \cap T_1^{-n}(U) \cap \cdots \cap T_N^{-n}(U) \ne \emptyset$). }
\end{definition}

Both types of transitvity are equivalent to their respective versions of dense hypercyclicity.  

\begin{proposition}{\normalfont{(\cite[Proposition 2.2]{BeP1}, \cite[Proposition 3.3]{BGJ})}} \label{prop:transitive}
    The operators $T_1, \dots, T_N \in B(X)$ with $N \ge 2$ are densely d-hypercyclic (respectively, densely s-hypercyclic) if and only if
    they are d-topologically transitive (respectively, s-topologically transitive).
\end{proposition}

There are multiple sufficient conditions to show a finite collection of operators are densely d-hypercyclic or densely s-hypercyclic.  The blow-up/collapse properties are one such method.

\begin{definition}\label{def:d-s-BlowUp}
    \normalfont{
    We say the operators $T_1, \dots, T_N \in B(X)$ with $N \ge 2$ satisfy the \textit{Strong Disjoint Blow-up/Collapse Property} (respectively, \textit{Disjoint Blow-up/Collapse Property}) provided for each integer $L \in \mathbb{N}$ (respectively, for the integer $L=1$ only) and for any non-empty open sets $W$, $U_{1-L}, \dots, U_0, U_1, \dots, U_N$ of $X$ with $0 \in W$, there exists an integer $n \in \mathbb{N}$ such that
    \begin{align}
        W \cap T_{1}^{-n}(U_1) \cap \cdots \cap T_N^{-n}(U_N) 
        \ne \emptyset
        \mbox{ and }
        U_{\ell} \cap T_{1}^{-n}(W) \cap \cdots \cap T_N^{-n}(W) 
        \ne \emptyset
        \mbox{ for integers $1-L \le \ell \le 0$.}
        \nonumber
    \end{align}
    If for any non-empty open sets $W, U,V$ in $X$ with $0 \in W$, there exists an integer $n \in \mathbb{N}$ such that
     \begin{align}
        W \cap T_{1}^{-n}(U) \cap \cdots \cap T_N^{-n}(U) 
        \ne \emptyset
        \mbox{ and }
        V \cap T_{1}^{-n}(W) \cap \cdots \cap T_N^{-n}(W) 
        \ne \emptyset,
        \nonumber
    \end{align}
    we say the operators $T_1, \dots, T_N$ satisfy the \textit{Simultaneous Blow-up/Collapse Property}}.
\end{definition}

Satisfying these blow-up/collapse properties imply the corresponding dense hypercyclicity.

\begin{proposition} {\normalfont{(\cite[Theorem 8]{Sa4},  \cite[Theorem 1.6]{BeP1}, \cite[Proposition 3.5]{BGJ})}}\label{prop:d-s-BlowUp} 
For operators $T_1, \dots, T_N \in B(X)$ with $N \ge 2$, the following assertions hold:
\begin{enumerate}
    \item[(i)] If the operators $T_1, \dots, T_N$ satisfy the Strong Disjoint Blow-up/Collapse Property, then they posses a dense d-hypercyclic linear manifold and so are densely d-hypercyclic.
    \item[(ii)] If the operators $T_1, \dots, T_N$ satisfy the Disjoint Blow-up/Collapse Property (respectively, Simultaneous Blow-Up/Collapse Property), then they are densely d-hypercyclic (repspectively, densely s-hypercyclic). 
\end{enumerate}
\end{proposition}

With Proposition \ref{prop:transitive} and the original definition of d-weakly mixing, we say operators $T_1, \dots, T_N \in B(X)$ with $N \ge 2$ are \textit{d-weakly mixing} (respectively, \textit{s-weakly mixing}) if the direct sums $T_1 \oplus T_1, \dots, T_N \oplus T_N$ are d-topologically transitive (respectively, s-topologically transitive).  The Disjoint Hypercyclicity Criterion is a sufficient condition for d-weakly mixing.  

\begin{definition} {\normalfont{\cite[Definition 2.5]{BeP1}}} \label{def:CoHC}
	\normalfont{The operators $T_1,\dots,T_N\in B(X)$ with $N \ge 2$ satisfy the \textit{Disjoint Hypercyclicity Criterion}, or \textit{d-Hypercyclicity Criterion} for short, with respect to the strictly increasing sequence $(n_k)$ of positive integers provided there exist dense subsets $X_0, X_1,\dots, X_N$ of $X$ and mappings
		$S_{j, k}:X_j \longrightarrow X$ for integers $k, j \in \mathbb{N}$ with $1\le j \le N$ satisfying the following conditions:
		\begin{enumerate}
		    \item[(a)]  for each $i$ with $1 \le i \le N$, we have $T_i^{n_k} \longrightarrow 0$ pointwise on $X_0$ as $k \longrightarrow \infty$,
		    \item[(b)] for each integer $j$ with $1 \le j \le N$, we have $S_{j,k} \longrightarrow 0$ pointwise on $X_j$ as $k \longrightarrow \infty$, and 
		    \item[(c)] for integers $i,j$ with $1 \le i,j \le N$, we have $(T_i^{n_k} S_{j,k} - \delta_{i,k} Id_{X_i}) \longrightarrow 0$ pointwise in $X_i$ as $k \longrightarrow \infty$.
		\end{enumerate}
		In general, we say the operators $T_1,\dots,T_N$ satisfy the d-Hypercyclicity Criterion if there exists some sequence $(n_k)$ for which the above conditions are satisfied.}
\end{definition}

The Disjoint Hypercyclicity Criterion not only ensures a finite collection of operators is d-weakly mixing, but it is equivalent to direct sums of those operators of arbitrary length are densely d-hypercyclic.

\begin{proposition}{\normalfont{\cite[Theorem 2.7]{BeP1}}}\label{prop:CoHCrit}
	For any operators $T_1,\dots, T_N \in B(X)$ with $N \ge 2$, the following assertions are equivalent: 
	\begin{enumerate}
	    \item[(i)] The operators $T_1, \dots, T_N$ satisfy the d-Hypercyclicity Criterion.
	    \item[(ii)] The operators $T_1, \dots, T_N$ are hereditarily densely d-hypercyclic for some strictly increasing sequence $(n_k)$ of positive integers.
	    \item[(iii)] For each integer $R \in \mathbb{N}$, the direct sum operators $\oplus_{r=1}^{R} T_1, \dots, \oplus_{r=1}^{R} T_N$ are d-topologically transitive.
	\end{enumerate}
	Furthermore, whenever the operators $T_1, \dots, T_N$ satisfy the d-Hypercyclicity Criterion, they posses a dense d-hypercyclic linear manifold and hence are densely d-hypercyclic.
\end{proposition}

The Simultaneous Hypercyclicity Criterion  serves a parallel role for s-hypercyclicity. For the necessary notation to state the criterion,
for any subset $A$ of a vector space, the convex hull conv$(A)$ is the smallest convex subset of of the vector space containing $A$.

\begin{definition}\cite[Definition 3.6]{BGJ} \label{def:sHC}
	\normalfont{
	The operators $T_1, \dots, T_N \in B(X)$ with $N \ge 2$ satisfy the \textit{Simultaneous Hypercyclicity Criterion}, or \textit{s-Hypercyclicity Criterion} for short, with respect to the strictly increasing sequence $(n_k)\subset\N$ of positive integers provided there exist a dense subset $X_0$ in $X$, subset $W_0$ in $\oplus_{i=1}^{N} X$ with 
	$\overline{W_0}\supset \Delta \left( \oplus_{i=1}^{N} X \right)$,	
    and mappings $R_k: W_0\longrightarrow X$ for each $k \in \mathbb{N}$ satisfying the following conditions:
    \begin{enumerate}
	    \item[(a)] for each integer $i$ with $1 \le i \le N$, we have. $T_i^{n_k}\longrightarrow 0$ pointwise on $X_0$ as $k \longrightarrow \infty$,
	    \item[(b)] the mappings $R_k \longrightarrow 0$ pointwise on $W_0$ as $k \longrightarrow \infty$, and
	    \item[(c)] for each vector $w=(w_1,\ldots,w_N)\in W_0\,$ and each integer $i$ with $1 \le i \le N$, there is a vector $y_i \in \mbox{conv}(\{w_1\ldots,w_N\})$ such that 
	    $T_i^{n_k}R_k w \longrightarrow y_j$ as $k \longrightarrow \infty$. 
    \end{enumerate}
    In general, we say the operators $T_1, \dots, T_N$ satisfy the s-Hypercyclicity Criterion if there exists some sequence $(n_k)$ for which the above conditions are satisfied.
}
\end{definition}

The following result about the s-Hypercyclicity Criterion parallels the d-Hypercyclicty Criterion result in Proposition \ref{prop:CoHCrit}.

\begin{proposition}{\normalfont{\cite[Theorem 3.7]{BGJ}}} \label{prop:sHC}
	For any operators $T_1,\dots,T_N\in B(X)$ with $N \ge 2$, the following assertions are equivalent:
	\begin{enumerate}
	    \item[(i)] The operators $T_1, \dots, T_N$ satisfy the s-Hypercyclicity Criterion.
	    \item[(ii)] The operators $T_1, \dots, T_N$ are hereditarily densely s-hypercyclic for some strictly increasing sequence $(n_k)$ of positive integers.
	    \item[(iii)] For each integer $R \in \mathbb{N}$, the direct sum operators $\oplus_{r=1}^{R} T_1, \dots, \oplus_{r=1}^{R} T_N$ are s-topologically transitive.
	\end{enumerate}
	In particular, the operators $T_1, \dots, T_N$ are densely s-hypercyclic whenever the satisfy the s-Hypercyclicity Criterion.
\end{proposition}

\section{Disjoint Hypercyclic Pseudo-Shifts}\label{DisjointPseudoShfits}

Pseudo-shift operators serve as a generalization of shifting operators on sequence spaces. 
In this section, we provide a characterization of the d-hypercyclic unilateral pseudo-shifts on $\ell^p(\mathbb{N})$, which unifies several of the existing characterizations of d-hypercyclic unilateral weight shifts under a common umbrella.  
In order to formally define unilateral pseudo-shifts on the Banach space $\ell^{p}(\mathbb{N})$ with $1 \le p < \infty$ over the real or complex  scalar field $\mathbb{K}$, consider the standard canonical basis $\{ e_m : m \in \mathbb{N} \}$ of $\ell^p(\mathbb{N})$.  Every vector $x = (\alpha_1, \alpha_2, \alpha_3, \dots)$ in $\ell^p(\mathbb{N})$ can be represented as the convergent series
$x = \sum_{m=1}^{\infty} \alpha_m e_m$.
Using the series representation of the a vector in $\ell^p(\mathbb{N})$, we define a unilateral pseudo-shift in the following manner.

\begin{definition}\label{def:pseudoshift}{\rm
	Let $f: \mathbb{N} \longrightarrow \mathbb{N}$ be a strictly increasing map with $f(1) > 1$ and let $\omega = (w_{f(m)})_{m \in \mathbb{N}}$ be a bounded, nonzero weight sequence.  For $1 \le p < \infty$, the \textit{unilateral pseudo-shift} $T_{f,\omega}$ on $\ell^p(\mathbb{N})$, induced by the map $f$ and weight sequence $\omega$, is the linear operator given by
	\begin{align}
	    T_{f, \omega} \left(
	    \sum_{m=1}^{\infty} \alpha_m e_m \right)
	    = \sum_{m=1}^{\infty} w_{f(m)} \alpha_{f(m)} e_m
	    \nonumber
	\end{align}
	for every vector $x = \sum_{m=1}^{\infty} \alpha_m e_m $ in $\ell^p(\mathbb{N})$.
}
\end{definition}

To assist with the examination of the dynamics of the unilateral pseudo-shift operators, for any integers $n,M \in \mathbb{N}$ and any subset $A \subset \N$, we define $[M] = \{1,2,\ldots,M\}$, $f(A) = \{f(m): m \in A\}$,  and $f^{n} = f \circ \ldots \circ f,$ which is $f$ composed with itself $n$ many times. Furthermore, let $f^{-n}$ denote the inverse of $f^{n}$ defined on its image $f^n(\mathbb{N})$.
Note that the $n$-th iterate $T_{f,\omega}^n$ of the unilateral pseudo-shift is given by
\begin{align}
     T_{f, \omega}^n \left(
	    \sum_{m=1}^{\infty} \alpha_m e_m \right)
	    = \sum_{m=1}^{\infty} \alpha_{f^n(m)} W_{m,n} e_m
	    \mbox{ \ where \ } W_{m,n} = \prod_{\nu=1}^{n} w_{f^{\nu}(m)}
	    \label{intro.3}
\end{align}
for every vector $x = \sum_{m=1}^{\infty} \alpha_m e_m$ in $\ell^p(\mathbb{N})$.

Before we state the d-hypercyclic characterization, we provide a sufficient condition for satisfying the Disjoint Blow-up/Collapse Property, which is strictly weaker condition for d-hypercyclicity than the Disjoint Hypercyclicity Criterion introduced in \cite{BeP1}. 

\begin{theorem}[Disjoint Blow-up/Collapse Criterion] \label{weakcriterion}
	For operators $T_1, \dots T_N \in B(X)$ with $N \ge 2$, suppose there exist a strictly increasing sequence $(n_k)$ of positive integers, a dense subset $X_0$ of $X$ and maps $S_k: \oplus_{i=1}^{N} X_0 \longrightarrow X$ which satisfy the following:
	\begin{enumerate}
	    \item[(i)] for each vector $x \in X_0$ and integer $i$ with $1 \le i \le N$, we have
	    $T_i^{n_k} x \longrightarrow 0$ as $k \longrightarrow \infty$,
	    \item[(ii)] for each $\epsilon > 0$, integer $K \in \mathbb{N}$ and vectors $x_1, \dots, x_N \in X_0$, there exists an integer $k \ge K$ such that
	    \begin{enumerate}
	        \item[(a)] $\| S_k(x_1, \dots, x_N) \| < \epsilon$ and
	        \item[(b)]  $\| T_i^{n_k} S_k(x_1, \dots, x_N) - x_i \| < \epsilon$, 
	        for integers $1 \le i \le N$.
	        \nonumber
	    \end{enumerate}
	\end{enumerate}
Then the operators $T_1, \dots, T_N$ satisfy the Strong Disjoint Blow-up/Collapse Property, and hence posses a dense d-hypercyclic manifold.
\end{theorem}

\begin{proof}
	To verify the operators $T_1, \dots, T_N$ satisfy the Strong Disjoint Blow-up/Collapse Property, let $L \in \mathbb{N}$ and let $W$, $U_{1-L}, \dots, U_{-1}, U_0, U_1, \dots U_N$ be any non-empty open sets in $X$ with $0 \in W$.  By the density of the set $X_0$, there exist vectors $x_{1-L}, \dots, x_{-1}, x_0, x_1, \dots, x_N$ in $X_0$ with $x_{\ell} \in U_{\ell} \cap X_0$ for integers $1-L \le \ell  \le N$.
	By assumptions (i) and (ii) and the fact $0 \in W$, we may select an integer $k \ge K$ such that
	\begin{align}
	    & T_i^{n_k} x_{\ell} \in W \mbox{ for integers $1\le i \le N$ and $1-L \le \ell \le 0$},
	    \nonumber \\
	    & S_k(x_1, \dots, x_N) \in W \mbox{ and }
	     T_i^{n_k}S_k(x_1, \dots, x_N) \in U_i \mbox{ for integers $1 \le i \le N$},
	    \nonumber 
	\end{align}
    Therefore, $x_{\ell} \in U_{\ell} \cap T_1^{-n_k}(W) \cap \cdots \cap T_N^{-n_k}(W)$ for each integer $\ell$ with $1-L \le \ell \le 0$ and $S_k(x_1, \dots, x_N) \in W \cap T_1^{-n_k}(U_1) \cap \cdots \cap T_N^{-n_k}(U_N)$.
\end{proof}

The Disjoint Blow-up/Collapse Criterion has the following equivalent formulation, which allows one to suppress the maps $S_k$.

\begin{remark}\label{AlternativeWeakcriterion}
    The operators $T_1, \dots T_N \in B(X)$ with $N \ge 2$ satisfy the Disjoint Blow-up/Collapse Criterion if and only if there is a dense set $X_0$ such that for each $\epsilon > 0$, integers $K, L \in \mathbb{N}$ and vectors $y_1, \dots, y_L$, $x_1, \dots, x_N \in X_0$, there exists an integer $n \ge K$ and a vector $z \in X$ with $\| z \| < \epsilon$ such that 
    \begin{align}
        \| T_i^n y_{\ell} \| < \epsilon \mbox{ \ and \ } \| T_i^n z - x_i \| < \epsilon \mbox{ for integers $1 \le i \le N$ and $1 \le \ell \le L$.}
        \nonumber
    \end{align}
\end{remark}

In \cite{BeMaSa}, B\`es et al.\!\! showed that finite collections of d-hypercyclic weighted shift operators never satisfy the Disjoint Hypercyclicity Criterion although they always satisfy the Disjoint Blow-up/Collapse Property. In the next Theorem, we show that  pseudo-shifts are d-hypercyclic if and only if they satisfy the Disjoint Blow-Up/Collapse Criterion.  The theorem also completely characterize d-hypercyclic unilateral pseudo-shifts in terms of their weight sequences. Since pseudo-shifts include weighted shifts as a sub-family and they never satisfy the Disjoint Hypercyclicity Criterion, this characterization shows that Disjoint Blow-Up/Collapse Criterion is strictly weaker than the Disjoint Hypercyclicity Criterion. 

\begin{theorem}\label{thm:characterisation}
    Let $1 \le p < \infty$ and $N \ge 2$.
    For each integer $i$ with $1 \le i \le N$, let $T_i = T_{f_i, \omega^{(i)}}$ be the unilateral pseudo-shift on $\ell^{p}(\mathbb{N})$ induced by the map $f_i: \mathbb{N} \longrightarrow \mathbb{N}$ and the weight sequence $\omega^{(i)} = (w_{f_i(m)}^{(i)})_{m \in \mathbb{N}}$ as defined in Definition \ref{def:pseudoshift}, and
    let
    $W_{m,n}^{(i)} = \prod_{\nu=1}^{n} w_{f_i^{\nu}(m)}^{(i)}$ as defined in (\ref{intro.3}).
    The following assertions are equivalent:
    \begin{enumerate}
        \item[(i)] The pseudo-shifts $T_1, \dots, T_N$ are d-hypercyclic.
        \item[(ii)] The pseudo-shifts $T_1, \dots, T_N$ satisfy the Disjoint Blow-up/Collapse Criterion.
        \item[(iii)] There exists a strictly increasing sequence $(n_k)$ of positive integers which satisfies the following:
        \begin{enumerate}
            \item[(a)] for each integer $i,m \in \mathbb{N}$ with $1 \le i \le N$, we have
            $\left| W_{m,n_k}^{(i)} \right| \ \longrightarrow \infty$ as $k \longrightarrow \infty$, and
            \item[(b)] for each $\epsilon > 0$, integers $K, M \in \mathbb{N}$, and  finite collection $\{ a_{j}^{(i)} : 1 \le j \le M, 1 \le i \le N \}$ of non-zero scalars in $\mathbb{K} \setminus \{ 0 \}$, there exists an integer $k \ge K$ such that for any integers $i,\ell$ with $1 \le i,\ell \le N$ and $i \ne \ell$, we have
            \begin{align}
                & \left| \frac{W_{f_{i}^{-n_k}(j),n_k}^{(i)}}{W_{f_{\ell}^{-n_k}(j),n_k}^{(\ell)}} \right| < \epsilon 
                \mbox{ \  whenever $j \in f_{\ell}^{n_k}([M]) \cap f_{i}^{n_k}(\mathbb{N}\setminus [M])$,}
                \nonumber
            \end{align}
            and
            \begin{align}
                & \left| \frac{W_{f_{i}^{-n_k}(j),n_k}^{(i)}}{W_{f_{\ell}^{-n_k}(j),n_k}^{(\ell)}} - \frac{a_{f_{i}^{-n_k}(j)}^{(i)}}{a_{f_{\ell}^{-n_k}(j)}^{(\ell)}}
                \right| < \epsilon 
                \mbox{ \ whenever $j \in f_{\ell}^{n_k}([M]) \cap f_{i}^{n_k}([M])$}.
                \nonumber
            \end{align}
        \end{enumerate}
    \end{enumerate}
\end{theorem}

Before we prove Theorem \ref{thm:characterisation}, we first establish a technical lemma, which will be utilized in several other results in Section \ref{DisjointPseudoShfits} and Section \ref{SimultaneousPseudoShifts}.

\begin{lemma}\label{lem:technical}
    Let $T_1, \dots, T_N$ with $N \ge 2$ be unilateral pseudo-shifts on $\ell^{p}(\mathbb{N})$ as defined in Theorem \ref{thm:characterisation}.  
    For integers $n, M \in \mathbb{N}$ and vectors $x_1, \dots, x_N$ given by
    $x_i = \sum_{m=1}^{M} a_{m}^{(i)} e_m$ with $a_1^{(i)}, \dots, a_{M}^{(i)} \ne 0$,
    there exists a vector $z \in \mbox{span} \{ e_m : m \in \mathbb{N} \}$ such that
    \begin{align}
         \| z \| & \le MN \Gamma \max \left\{ \left| W_{m,n}^{(\ell)} \right|^{-1} : 1 \le \ell \le N, 1 \le m \le M \right\},
        \label{lem:technical.a} \\
        \| T_1^n z_n - x_1 \| 
        & \le M N \Gamma \max \left\{ \left| \frac{W_{f_1^{-n}(j),n}^{(1)}}{W_{f_{\ell}^{-n}(j),n}^{(\ell)}} \right| : 2 \le \ell \le N, j \in f_{\ell}^n([M]) \cap 
        f_1^n(\mathbb{N}\setminus [M])\right\},
        \label{lem:technical.b}
    \end{align}
    and for each integer $i$ with $2 \le i \le N$, we have
    \begin{align}
        \| T_i^n z - x_i \| 
        & \le M N \Gamma \max \left\{ \left| \frac{W_{f_i^{-n}(j),n}^{(i)}}{W_{f_{\ell}^{-n}(j),n}^{(\ell)}} 
        - \frac{a_{f_i^{-n}(j)}^{(i)}}{a_{f_{\ell}^{-n}(j)}^{(\ell)}}
        \right| : 1 \le \ell \le i-1, j \in f_{\ell}^n([M])\cap f_i^n([M])\right\} 
        \label{lem:technical.c} \\
        & \ \ \ \ + M N \Gamma
        \max \left\{ 
        \left| \frac{W_{f_i^{-n}(j),n}^{(i)}}{W_{f_{\ell}^{-n}(j),n}^{(\ell)}} \right| : \ell \ne i, j \in f_{\ell}^n([M]) \cap 
        f_i^{n}(\mathbb{N}\setminus [M])
        \right\}.
        \nonumber
    \end{align}
    where $\Gamma = \max \{ \| x_i \| : 1 \le i \le N \}$ and $W_{m,n}^{(i)} = \prod_{\nu=1}^{n} w_{f_i^{\nu}(m)}^{(i)}$.
\end{lemma}

\begin{proof}
Consider the vector $z \in \mbox{span} \{ e_m : m \in \mathbb{N} \}$ given by
	\begin{align}
	    z = \sum_{\ell=1}^{N} \left( 
	    \sum_{j \in A_{\ell}} \frac{a_{f_{\ell}^{-n}(j)}^{(\ell)}}{W_{f_{\ell}^{-1}(j),n}^{(\ell)}} \ e_j
	    \right).
	    \nonumber
	\end{align}
	where $A_{1}, \dots, A_{N}$ are the disjoint sets given by
	\begin{align}
        A_{\ell} = 
        \begin{cases}
            f_1^n([M]), & \text{ if $\ell=1$} \\
            f_{\ell}^n([M]) \setminus \bigcup_{i=1}^{\ell-1} f_i^n([M]), 
            & \text{ if $2 \le \ell \le N$}
        \end{cases}.
        \nonumber
    \end{align}
	
	To verify we have the desired vector $z$, observe that by the definition of constant $\Gamma$ and the fact that $A_{\ell} \subseteq f_{\ell}^{n}([M])$, it follows that
	\begin{align}
        \| z \|
        & \le \sum_{\ell=1}^{N}  
        \sum_{j \in A_{\ell}} \left|
        \frac{a_{f_{\ell}^{-n}(j)}^{(\ell)}}{W_{f_{\ell}^{-1}(j),n}^{(\ell)}} \right| 
        \nonumber \\
        & \le M N \Gamma 
        \max \left\{ \left| W_{f_{\ell}^{-n}(j),n}^{(\ell)} \right|^{-1} : 1 \le \ell \le N, j \in A_{\ell} \right\}
        \nonumber \\
        & \le M N \Gamma 
        \max \left\{ \left| W_{m,n}^{(\ell)} \right|^{-1} : 1 \le \ell \le N, 1 \le m \le M \right\}.
        \nonumber
    \end{align}
    Next, observe that $f_1^n([M]) = A_1$ and $A_{\ell} \cap f_1^n([M]) = \emptyset$ for each integer $\ell$ with $2 \le \ell \le N$, and so we can write
    \begin{align}
        T_1^n z 
        & = \sum_{\ell=1}^{N} \left( 
	    \sum_{j \in A_{\ell}} \frac{a_{f_{\ell}^{-n}(j)}^{(\ell)}}{W_{f_{\ell}^{-1}(j),n}^{(\ell)}} T_1^n e_j
	    \right) 
	    \nonumber \\
	    & = \sum_{j \in A_{1}} a_{f_1^{-n}(j)}^{(1)} e_{f_1^{-n}(j)}
	    + \sum_{\ell=2}^{N} \left( \sum_{j \in A_{\ell}} \frac{a_{f_{\ell}^{-n}(j)}^{(\ell)} W_{f_1^{-n}(j),n}^{(1)}}{W_{f_{\ell}^{-n}(j),n}^{(\ell)}} e_{f_1^{-n}(j)}
	    \right)
	    \nonumber \\
	    & = x_1 + 
	    \sum_{\ell=2}^{N} \left( \sum_{j \in A_{\ell} \cap f_1^n(\mathbb{N} \setminus [M])} \frac{a_{f_{\ell}^{-n}(j)}^{(\ell)} W_{f_1^{-n}(j),n}^{(1)}}{W_{f_{\ell}^{-n}(j),n}^{(\ell)}} e_{f_1^{-n}(j)}
	    \right).
	    \nonumber
    \end{align}
    It follows by the definition of the constant $\Gamma$ and $A_{\ell} \subseteq f_{\ell}^n([M])$ that
    \begin{align}
        \| T_1^n z - x_1 \|
        & \le  
        \sum_{\ell = 2}^{N} \left( \sum_{j \in A_{\ell} \cap f_1^n(\mathbb{N} \setminus [M])} \left| \frac{a_{f_{\ell}^{-n}(j)}^{(\ell)} W_{f_1^{-n}(j),n}^{(1)}}{W_{f_{\ell}^{-n}(j),n}^{(\ell)}}
	    \right| \right)
	    \nonumber \\
	   & \le M N \Gamma 
	   \max \left\{ \left| \frac{W_{f_1^{-n}(j),n}^{(1)}}{W_{f_{\ell}^{-n}(j),n}^{(\ell)}} \right| : 2 \le \ell \le N, j \in A_{\ell} \cap f_1^n(\mathbb{N} \setminus [M])\right\}
	   \nonumber \\
	   & \le M M \Gamma 
	   \max \left\{ \left| \frac{W_{f_1^{-n}(j),n}^{(1)}}{W_{f_{\ell}^{-n}(j),n}^{(\ell)}} \right| : 2 \le \ell \le N, j \in f_{\ell}^n([M]) \cap f_1^n(\mathbb{N}\setminus [M])\right\}.
	   \nonumber
    \end{align}
    Lastly, for each integer $i$ with $2 \le i \le N$, note that $A_{\ell} \cap f_i^n([M]) = \emptyset$ for each integer $\ell$ with $i+1 \le \ell \le N$, and $f_i^n([M]) = \bigcup_{\ell=1}^{i} A_{\ell} \cap f_i^n([M])$ with $A_i \cap f_i^{n}([M]) = A_i$.  Thus, we can write
    \begin{align}
        x_i & = 
        \sum_{j \in A_i} a_{f_i^{-n}(j)}^{(i)} e_{f_i^{-n}(j)}
        + \sum_{\ell=1}^{i-1} \left(
        \sum_{j \in A_{\ell} \cap f_i^n([M])} a_{f_i^{-n}(j)}^{(i)} e_{f_i^{-n}(j)}
        \right),
        \label{lem:technical.10}
    \end{align}
    and
    \begin{align}
        T_i^n z 
        & = \sum_{\ell=1}^{N} \left( 
	    \sum_{j \in A_{\ell}} \frac{a_{f_{\ell}^{-n}(j)}^{(\ell)}}{W_{f_{\ell}^{-1}(j),n}^{(\ell)}} T_i^n e_j
	    \right) 
	    \label{lem:technical.11} \\
	    & = \sum_{j \in A_i} a_{f_i^{-n}(j)}^{(i)} e_{f_i^{-n}(j)}
	    + \sum_{\ell \ne i} \left( \sum_{j \in A_{\ell}}
	    \frac{a_{f_{\ell}^{-n}(j)}^{(\ell)} W_{f_i^{-n}(j),n}^{(i)}}{W_{f_{\ell}^{-n}(j),n}^{(\ell)}} e_{f_i^{-n}(j)}
	    \right).
	    \nonumber 
	    \\
	    & = \sum_{j \in A_i} a_{f_i^{-n}(j)}^{(i)} e_{f_i^{-n}(j)}
	    + \sum_{\ell=1}^{i-1} \left(
	    \sum_{j \in A_{\ell} \cap f_i^n([M])}
	    \frac{a_{f_{\ell}^{-n}(j)}^{(\ell)} W_{f_i^{-n}(j),n}^{(i)}}{W_{f_{\ell}^{-n}(j),n}^{(\ell)}} e_{f_i^{-n}(j)}
	    \right)
	    \nonumber 
	    \\
	    & \quad \ \ \ \ + 
	    \sum_{\ell \ne i} \left(
	    \sum_{j \in A_{\ell} \cap f_i^n(\mathbb{N} \setminus [M])}
	    \frac{a_{f_{\ell}^{-n}(j)}^{(\ell)} W_{f_i^{-n}(j),n}^{(i)}}{W_{f_{\ell}^{-n}(j),n}^{(\ell)}} e_{f_i^{-n}(j)}
	    \right).
	    \nonumber
    \end{align}
    Combining equations (\ref{lem:technical.10}) and (\ref{lem:technical.11}) with the definition of $\Gamma$ and the fact that $A_{\ell} \subseteq f_{\ell}^n([M])$ gives us
    \begin{align}
        \| T_i^n z - x_i \|
        & \le
        \sum_{\ell=1}^{i-1} \sum_{j \in A_{\ell} \cap f_i^n([M])} 
        \left| a_{f_{\ell}^{-n}(j)}^{(\ell)} \right| \left|
        \frac{ W_{f_i^{-n}(j),n}^{(i)}}{W_{f_{\ell}^{-n}(j),n}^{(\ell)}} e_{f_i^{-n}(j)}
        - \frac{a_{f_i^{-n}(j)}^{(i)}}{a_{f_{\ell}^{-n}(j)}^{(\ell)}}
        \right|
        \nonumber \\
        & \quad \ \ \ \ +
        \sum_{\ell \ne i} \sum_{j \in A_{\ell} \cap f_i^n(\mathbb{N} \setminus [M])}
        \left|
        \frac{a_{f_{\ell}^{-n}(j)}^{(\ell)} W_{f_i^{-n}(j),n}^{(i)}}{W_{f_{\ell}^{-n}(j),n}^{(\ell)}}
        \right|
        \nonumber \\
        & \le M N \Gamma \max \left\{ \left| \frac{W_{f_i^{-n}(j),n}^{(i)}}{W_{f_{\ell}^{-n}(j),n}^{(\ell)}} 
        - \frac{a_{f_i^{-n}(j)}^{(i)}}{a_{f_{\ell}^{-n}(j)}^{(\ell)}}
        \right| : 1 \le \ell \le i-1, j \in A_{\ell} \cap f_i^n([M])\right\} 
        \nonumber \\
        & \ \ \ \ + M N \Gamma
        \max \left\{ 
        \left| \frac{W_{f_i^{-n}(j),n}^{(i)}}{W_{f_{\ell}^{-n}(j),n}^{(\ell)}} \right| : \ell \ne i, j \in A_{\ell} \cap f_i^{n}(\mathbb{N} \setminus [M])
        \right\}
        \nonumber \\
        & \le M N \Gamma \max \left\{ \left| \frac{W_{f_i^{-n}(j),n}^{(i)}}{W_{f_{\ell}^{-n}(j),n}^{(\ell)}} 
        - \frac{a_{f_i^{-n}(j)}^{(i)}}{a_{f_{\ell}^{-n}(j)}^{(\ell)}}
        \right| : 1 \le \ell \le i-1, j \in f_{\ell}^n([M])\cap f_i^n([M])\right\} 
        \nonumber \\
        & \ \ \ \ + M N \Gamma
        \max \left\{ 
        \left| \frac{W_{f_i^{-n}(j),n}^{(i)}}{W_{f_{\ell}^{-n}(j),n}^{(\ell)}} \right| : \ell \ne i, j \in f_{\ell}^n([M]) \cap f_i^{n}(\mathbb{N} \setminus [M])
        \right\},
        \nonumber
    \end{align}
    which concludes the proof of the lemma.
\end{proof}

Now we are ready to prove Theorem \ref{thm:characterisation}.

\begin{proof}
	To begin, we establish (iii) implies (ii) by assuming there is a strictly increasing sequence $(n_k)$ of positive integers which satisfy conditions (a) and (b) in statement (iii), and we verify the unilateral pseudo-shifts $T_1, \dots, T_N$ satisfy the Disjoint Blow-up/Collapse Criterion by establishing the equivalent condition given in Remark \ref{AlternativeWeakcriterion}.  To this end, let $X_0 = \mbox{span} \{ e_m : m \in \mathbb{N} \}$, which is dense in $\ell^p(\mathbb{N})$.
	Consider an arbitrary $\epsilon > 0$, integers $K,L \in \mathbb{N}$ and vectors $y_1, \dots, y_L$, $x_1, \dots, x_N \in X_0$.  Note we can write each vector $x_i$ as
	\begin{align}
	    x_i = \sum_{m=1}^{M} a_{m}^{(i)} e_m,
	    \mbox{ \ where $M = \max \left\{ \mbox{deg}(y_{\ell}), \mbox{deg}(x_i) : 1 \le \ell \le L, 1 \le i \le N \right\}$.}
	    \label{thm:charact.2}
	\end{align}
	Without loss of generality, we may assume
	$a_{1}^{(i)}, \dots,  a_{M}^{(i)}\ne 0$ for integers $1 \le i \le N$.  As a verification, note that one may select vectors $\widetilde{x}_1, \dots, \widetilde{x}_N$ with $\widetilde{x}_1 = \sum_{m=1}^{M} \tilde{a}_{m}^{(i)} e_m$ such such that
    $\| y_i - x_i \| < \frac{\epsilon}{2}$  and $\tilde{a}_{1}^{(i)}, \dots, \tilde{a}_{M}^{(i)} \ne 0$. 
	Then whenever $\| T_i^n z - \widetilde{x}_i \| < \frac{\epsilon}{2}$, it follows that
	$\| T_i^n z - x_i \| \le \| T_i^n z - \widetilde{x}_i \| + \| \widetilde{x}_i - x_i \| < \epsilon$.
	
	Returning to the proof, by our assumptions, we may select in integer $n = n_k$ with $k \ge \max \{ K, M \}$ such that
	\begin{align}
	    \left| W_{m,n}^{(i)} \right| & > \frac{M N \Gamma}{\epsilon},
	    \mbox{ \ for integers $1 \le i \le N$ and $1 \le m \le M$},
	    \label{thm:charact.6} \\
	    \left| \frac{W_{f_i^{-n}(j),n}^{(i)}}{W_{f_{\ell}^{-n}(j),n}^{(\ell)}} \right| 
	    & < \frac{\epsilon}{2 M N \Gamma}, \mbox{ \ whenever $i \ne \ell$
	    and $j \in f_{\ell}^n([M]) \cap f_i^n(\mathbb{N} \setminus [M])$},
	    \label{thm:charact.7} \\
	    \left| \frac{W_{f_i^{-n}(j),n}^{(i)}}{W_{f_{\ell}^{-n}(j),n}^{(\ell)}}  - \frac{a_{f_i^{-n}(j)}^{(i)}}{a_{f_{\ell}^{-n}(j)}^{(\ell)}} \right|  
	    & < \frac{\epsilon}{2 M N \Gamma}, \mbox{ \ whenever $i \ne \ell$
	    and $j \in f_{\ell}^n([M]) \cap f_i^n([M])$},
	    \label{thm:charact.8}
	\end{align}
	where $\Gamma = \max \{ \| x_i \| : 1 \le i \le N \}$.
	By Lemma \ref{lem:technical} using the integers $M$, $n=n_k$ and vectors $x_1, \dots, x_N$, there is a vector $z \in X_0$ which satisfies inequalities (\ref{lem:technical.a}), (\ref{lem:technical.b}) and (\ref{lem:technical.c}).
	
	To verify we have the desired integer $n = n_k$ and vector $z$ required in Remark \ref{AlternativeWeakcriterion}, observe that since the strictly increasing map $f_i$ satisfies $f_i(1) > 1$, we have 
	\begin{align}
        \mbox{span} \{ e_m : 1 \le m \le n\} \subseteq \mbox{span} \{ e_m : j \in \mathbb{N} \setminus f_i^n(\mathbb{N}) \} = \mbox{Ker}(T_i^n).
        \label{thm:charact.11.1}
    \end{align}
    By the definition $M$ in (\ref{thm:charact.2}) and the fact $n = n_k \ge k \ge \max\{ K, M \}$, it immediately follows that $\| T_i^n y_{\ell} \| = 0 < \epsilon$ for any integers $i,\ell$ with $1 \le i \le N$ and $1 \le \ell \le L$.   
    Next, note that by inequality (\ref{lem:technical.a}) in Lemma \ref{lem:technical} and inequality (\ref{thm:charact.6}), we get
    \begin{align}
        \| z \| 
        & \le MN \Gamma \max \left\{ \left| W_{m,n}^{(\ell)} \right|^{-1} : 1 \le \ell \le N, 1 \le m \le M \right\} < \epsilon,
        \nonumber
    \end{align}
    and by inequality (\ref{lem:technical.b}) in Lemma \ref{lem:technical} and inequality (\ref{thm:charact.7}),
    \begin{align}
        \| T_1^n z - x_1 \|
        & \le M N \Gamma \max \left\{ \left| \frac{W_{f_1^{-n}(j),n}^{(1)}}{W_{f_{\ell}^{-n}(j),n}^{(\ell)}} \right| : 2 \le \ell \le N, j \in f_{\ell}^n([M]) \cap f_1^n((\mathbb{N} \setminus [M])\right\}
        < \epsilon.
        \nonumber
    \end{align}
    Likewise, for each integer $i$ with $2 \le i \le N$, we have
    \begin{align}
        \| T_i^n z - x_i \|
        & \le M N \Gamma \max \left\{ \left| \frac{W_{f_i^{-n}(j),n}^{(i)}}{W_{f_{\ell}^{-n}(j),n}^{(\ell)}} 
        - \frac{a_{f_i^{-n}(j)}^{(i)}}{a_{f_{\ell}^{-n}(j)}^{(\ell)}}
        \right| : 1 \le \ell \le i-1, j \in f_{\ell}^n([M])\cap f_i^n([M])\right\} 
        \nonumber  \\
        & \ \ \ \ + M N \Gamma
        \max \left\{ 
        \left| \frac{W_{f_i^{-n}(j),n}^{(i)}}{W_{f_{\ell}^{-n}(j),n}^{(\ell)}} \right| : \ell \ne i, j \in f_{\ell}^n([M]) \cap f_i^{n}( \mathbb{N} \setminus [M]) 
        \right\}, \mbox{ by (\ref{lem:technical.c}) in Lemma \ref{lem:technical}} 
        \nonumber \\
        & < \frac{\epsilon}{2} + \frac{\epsilon}{2}, 
        \mbox{ by (\ref{thm:charact.7}) and (\ref{thm:charact.8})}
        \nonumber \\
        & = \epsilon,
        \nonumber 
    \end{align}
    which concludes the proof of (iii) implies (ii).

    Since (ii) implies (i) follows directly from the Disjoint Blow-up/Collapse Criterion, we finish the proof of Theorem \ref{thm:characterisation} by establishing (i) implies (iii).  Suppose $x = \sum_{m=1}^{\infty} \alpha_m e_m$ in $\ell^p(\mathbb{N})$ is a disjoint hypercyclic vector for the unilateral pseudo-shifts $T_1, \dots, T_N$.  Note that by (\ref{thm:charact.11.1}) we have  $x \not\in \mbox{span} \{ e_m : m \in \mathbb{N} \}$, and so we can find a sequence $(\rho_k)$ of positive scalars such that
    \begin{align}
        \rho_k \longrightarrow 0 \mbox{ \ \ and \ \ } 
        \frac{\rho_k}{\sup \{ |\alpha_m| : m \ge k \} } \longrightarrow \infty \mbox{ \ as $k \longrightarrow \infty$.}
        \label{thm:charact.28}
    \end{align}
    Since $x$ is disjoint hypercyclic vector for the unilateral pseudo-shifts $T_1, \dots, T_N$, there exists a strictly increasing sequence $(n_k)$ of positive integers for which
    \begin{align}
        & \{ (T_1^{n_k} x, \dots, T_N^{n_k}x ) : k \ge 1 \} \mbox{ is dense in $\oplus_{i=1}^{N} \ell^p(\mathbb{N})$, and }
        \label{thm:charact.29} \\
        & \min \left\{ |\langle T_i^{n_k}x, e_m \rangle |: 1 \le i \le N, 1 \le m \le k \right\} > \rho_k \mbox{ for integers $k \ge 1$,}
        \label{thm:charact.30}
    \end{align}
    where $\langle \ \cdot \ , e_m \rangle$ is the linear functional on $\ell^p(\mathbb{N})$ associated with the vector $e_m$ in the dual space $\ell^p(\mathbb{N})^* = \ell^q(\mathbb{N})$.
    For any integers $i, m \in \mathbb{N}$ with $1 \le i \le N$, by the series representation of $T_i^{n_k} x$ given in (\ref{intro.3}) and by (\ref{thm:charact.30}), we have
    \begin{align}
        \rho_k <  |\langle T_i^{n_k}x, e_m \rangle |
        = | \alpha_{f_i^n(m)} | \left| W_{m,n_k}^{(i)} \right| 
        \le \left| W_{m,n_k}^{(i)} \right| \sup \{ | \alpha_m | : m \ge k \} \mbox{ \ for integers $k \ge m$.}
        \label{thm:charact.31}
    \end{align}
    From (\ref{thm:charact.28}), it follows that $\left| W_{m,n_k}^{(i)} \right| \longrightarrow \infty$ as $k \longrightarrow \infty$ as desired.  
    
    Next, consider $\epsilon > 0$, integers $K,M \in \mathbb{N}$ and finite collection $\{ a_{j}^{(i)} : 1 \le i \le N, 1 \le j \le M \}$ of non-zero scalars in $\mathbb{K} \setminus \{ 0 \}$.  Select a constant $C > 0$ such that
    \begin{align}
        1 < C \min \left\{ \left| a_{j}^{(i)} \right| : 1 \le i \le N, 1 \le j \le M \right\} 
        \mbox{ \ and \ }
        \max \left\{ \left| \frac{  a_{j}^{(i)} }{  a_{m}^{(\ell)}} \right| : 1 \le i,\ell \le N, 1 \le j,m \le M \right\} < C,
        \label{thm:charact.32}
    \end{align}
    and select a positive constant $\eta$  with  
    $0 < \frac{(1 + C) \eta}{1 - \eta} < \epsilon$.  By the density of the set in (\ref{thm:charact.29}), there exists an integer $k \ge K$ such that $\| T_i^{n_k}x - \sum_{m=1}^{M} C a_{m}^{(i)} e_m \| < \eta$ for each integer $i$ with $1 \le i \le N$.  Again by the series representation of $T_i^{n_k} x$ in (\ref{intro.3}), we have
    \begin{align}
        \left| \alpha_j W_{f_i^{-n_k}(j), n_k}^{(i)} - C a_{f_i^{-n_k}(j)}^{(i)} \right| 
        \le \left\| T_i^{n_k}x - \sum_{m=1}^{M} C a_{m}^{(i)} e_m \right\| 
        < \eta
        \mbox{ whenever $j \in f_i^{n_k}([M])$}
        \label{thm:charact.33}.
    \end{align}
    From the above inequality, the selection of the constant $C$ ensures $\alpha_j \ne 0$ and 
    $1 - \eta < \left| \alpha_j W_{f_i^{-n_k}(j), n_k}^{(i)}  \right|$ 
    for any integer $j \in f_i^{n_k}([M])$.  Also, note the inequality $\| T_i^{n_k}x - \sum_{m=1}^{M} C a_{m}^{(i)} e_m \| < \eta$ implies that
     \begin{align}
        \left| \alpha_j W_{f_i^{-n_k}(j), n_k}^{(i)} \right| 
        \le \left\| T_i^{n_k}x - \sum_{m=1}^{M} C a_{m}^{(i)} e_m \right\| 
        < \eta
        \mbox{ whenever $j \in f_i^{n_k}(\mathbb{N} \setminus [M])$}
        \label{thm:charact.34}.
    \end{align}
    For any integer $j \in f_{\ell}^{n_k}([M]) \cap f_i^{n_k}(\mathbb{N} \setminus [M])$ with $\ell \ne i$, by  (\ref{thm:charact.33}) and (\ref{thm:charact.34}), we have
    \begin{align}
        \left| \frac{W_{f_i^{-n_k}(j), n_k}^{(i)}}{W_{f_{\ell}^{-n_k}(j), n_k}^{(\ell)}} \right|
        = \frac{ \left| \alpha_j W_{f_i^{-n_k}(j), n_k}^{(i)} \right|}{ \left| \alpha_j W_{f_{\ell}^{-n_k}(j), n_k}^{(\ell)} \right|} 
        < \frac{\eta}{1-\eta}
        < \epsilon.
        \label{thm:charact.35}
    \end{align}
    For any integer $j \in f_{\ell}^{n_k}([M]) \cap f_i^{n_k}([M])$ with $\ell \ne i$, observe that
    \begin{align}
        \left| \frac{W_{f_i^{-n_k}(j), n_k}^{(i)}}{W_{f_{\ell}^{-n_k}(j), n_k}^{(\ell)}} - \frac{a_{f_i^{-n_k}(j)}^{(i)}}{a_{f_{\ell}^{-n_k}(j)}^{(\ell)}} \right|
        & \le 
        \frac{1}{ |\alpha_j W_{f_{\ell}^{-n_k}(j), n_k}^{(\ell)} | }
        \left| \alpha_j W_{f_i^{-n_k}(j), n_k}^{(i)} - 
        \frac{a_{f_i^{-n_k}(j)}^{(i)}}{a_{f_{\ell}^{-n_k}(j)}^{(\ell)}}
        \alpha_j W_{f_{\ell}^{-n_k}(j), n_k}^{(\ell)}
        \right|
        \label{thm:charact.36} \\
         & \le
        \frac{1}{ |\alpha_j W_{f_{\ell}^{-n_k}(j), n_k}^{(\ell)} | }
        \left|\alpha_j W_{f_i^{-n_k}(j), n_k}^{(i)}  - C a_{f_i^{-n_k}(j)}^{(i)} \right| 
        \nonumber \\
        & \ \ \ \ + 
        \frac{1}{ |\alpha_j W_{f_{\ell}^{-n_k}(j), n_k}^{(\ell)} | }
        \left| \frac{a_{f_i^{-n_k}(j)}^{(i)}}{a_{f_{\ell}^{-n_k}(j)}^{(\ell)}} \right|
        \left|C a_{f_{\ell}^{-n_k}(j)}^{(\ell)} -
        \alpha_j W_{f_{\ell}^{-n_k}(j), n_k}^{(\ell)}
        \right|
        \nonumber \\
         & < \frac{\eta}{1-\eta} + \frac{C \eta}{1-\eta}, \mbox{ by (\ref{thm:charact.32}), (\ref{thm:charact.33}) and (\ref{thm:charact.34})}
         \nonumber \\
         & < \epsilon.
         \nonumber 
    \end{align}
    
\end{proof}

Since satisfying the Disjoint Blow-up/Collapse Criterion for a finite collection of operators implies the existence of a dense d-hypercyclic manifold, it follows from Theorem \ref{thm:characterisation} that d-hypercyclicity and dense d-hypercyclicity are also equivalent concepts for unilateral pseudo-shifts.

\begin{corollary}\label{cor:densemanifold}
Let $T_1, \dots, T_N$ with $N \ge 2$ be unilateral pseudo-shifts on $\ell^{p}(\mathbb{N})$ as defined in Theorem \ref{thm:characterisation}.  The pseudo-shifts $T_1, \dots, T_N$ are d-hypercyclic if and only if they possess a dense d-hypercyclic manifold.
\end{corollary}

Shkarin, et al.\!\! \cite[Theorem 3.4]{SaSh} provided examples of d-hypercyclic operators whose set of d-hypercyclic vectors is contained within a finite dimensional subspace, and so d-hypercyclicity and dense d-hypercyclicity are not necessarily equivalent conditions.

A unilateral weighted backward shift is an example of a unilateral pseudo-shift with the strictly increasing map $f:\mathbb{N} \longrightarrow \mathbb{N}$ given by $f(m) = m+1$.  Utilizing this observation, the characterization of d-hypercyclic unilateral weighted backward shifts established by B\`es, et al.\! \cite{BeMaSa} is now a corollary of Theorem \ref{thm:characterisation}. 

\begin{corollary}\cite[Theorem 2.2]{BeMaSa} \label{crl:weightedshifts}
Let $1 \le p < \infty$ and let $N \ge 2$.  For each integer $i$ with $1 \le i \le N$, let $T_i$ be a unilateral weighted backward shift on $\ell^p(\mathbb{N})$ with the weight sequence $\omega^{(i)} = (w_m^{(i)})_{m \in \mathbb{N}}$. For integers $i,m,n \in \mathbb{N}$ with $2 \le i \le N$, define
\begin{align}
    \alpha_{m,n}^{(i)} = \prod_{\nu = 1}^{n} 
    \frac{w_{\nu + m}^{(i)}}{w_{\nu + m}^{(1)}}.
    \nonumber 
\end{align}
The following assertions are equivalent:
\begin{enumerate}
	\item[(i)] The shifts $T_1, \ldots, T_N$ are $d$-hypercyclic;
	\item[(ii)] The shifts $T_1, \ldots, T_N$ satisfy the Disjoint Blow up/Collapse Property;
	\item[(iii)] There exists a strictly increasing sequence $(n_k)$ of positive integers such that for each integer $m \in \mathbb{N}$, we have
	\begin{align} 
		\left|\prod_{\nu = 1}^{n_k} w_{\nu+ m}^{(1)}\right| \longrightarrow \infty \mbox{\it{ as }} k \longrightarrow \infty,
		\nonumber
	\end{align} 
	and the set 
	\begin{align}
		\left\{ (\alpha^{(2)}_{1,n_k},\alpha^{(3)}_{1,n_k},\ldots,\alpha^{(N)}_{1,n_k},\alpha^{(2)}_{2,n_k},
		\alpha^{(3)}_{2,n_k},\ldots,\alpha^{(N)}_{2,n_k},\ldots):k\ge 1 \right\}
		\nonumber
	\end{align} 
	is dense in $\mathbb{K}^{\mathbb{N}}$ with respect to the product topology.
\end{enumerate}
\end{corollary}

\begin{proof}
    As noted before the statement of Corollary \ref{crl:weightedshifts}, for each integer $i$ with $1 \le i \le N$, the unilateral weighted backward shift $T_i$ with weight sequence $\omega^{(i)} = (w_m^{(i)})_{m \in \mathbb{N}}$ can be expressed as the unilateral pseudo-shift $T_{f_i,\omega^{(i)}}$ where $f_i(n) = n+1$ for each $n \in \mathbb{N}$.  It clearly follows from Theorem \ref{thm:characterisation} and the Disjoint Blow-Up/Collapse Criterion that the shifts $T_1, \dots, T_N$ are d-hypercyclic if and only if they satisfy the Disjoint Blow-Up/Collapse Property.  Thus, it remains to show statement (iii) in Theorem \ref{thm:characterisation} is equivalent to statement (iii) in Corollary \ref{crl:weightedshifts}.  To this end, observe that one can easily show $W_{m,n}^{(i)} = \prod_{\nu = 1}^{n} w_{f_i^{\nu}(m)}^{(i)} = \prod_{\nu=1}^{n} w_{\nu + m}^{(i)}$ and 
    \begin{align}
        \frac{W_{m,n}^{(i)}}{W_{m,n}^{(\ell)}} =
        \begin{cases}
            \frac{1}{\alpha_{m,n}^{(\ell)}}, 
            & \text{ if $i=1$ and $2 \le \ell \le N$} \\
            \frac{\alpha_{m,n}^{(i)}}{\alpha_{m,n}^{(\ell)}},
            & \text{ if $2 \le i,\ell \le N$ and $i \ne \ell$} \\
            \alpha_{m,n}^{(i)},
            & \text{ if $2 \le i \le N$ and $\ell=1$}
        \end{cases}.
        \nonumber 
    \end{align}
    Furthermore, we have $f_{\ell}^n([M])\cap f_i^n(\mathbb{N}\setminus [M]) = \emptyset$ and $f_{\ell}^n([M]) \cap f_i^n([M]) = \{ n+1, \dots, n+M \}$ for any integers $i,\ell, n,M \in \mathbb{N}$ with $1 \le i,\ell \le N$.  Thus, the density of the set 
    \begin{align}
        \left\{ (\alpha^{(2)}_{1,n_k},\alpha^{(3)}_{1,n_k},\ldots,\alpha^{(N)}_{1,n_k},\alpha^{(2)}_{2,n_k},\alpha^{(3)}_{2,n_k}, \ldots,\alpha^{(N)}_{2,n_k},\ldots):k\ge 1
        \right\}
        \label{weightedshifts.2}
	\end{align}
	in $\mathbb{K}^{\mathbb{N}}$ with respect the strictly increasing sequence $(n_k)$ is equivalent to the density condition (b) in statement (iii) of Theorem \ref{thm:characterisation} with regard to the same sequence $(n_k)$.  Moreover, note that 
	$W_{m,n}^{(i)} = \alpha_{m,n}^{(i)} W_{m,n}^{(1)}$ 
	for each integer $i$ with $2 \le i \le N$.  Thus, the density of the set in (\ref{weightedshifts.2}) together with 
	$|W_{m,n_k}^{(1)}| \longrightarrow \infty$ as $k \longrightarrow \infty$ implies there is a subsequence $(n_{k_{r}})$ for which $|W_{m,n_{k_r}}^{(i)}| \longrightarrow \infty$ as $k \longrightarrow \infty$ for each integer $i$ with $1 \le i \le N$ and the set
	$\{(\alpha^{(2)}_{1,n_{k_r}},\alpha^{(3)}_{1,n_{k_r}},\ldots,\alpha^{(N)}_{1,n_{k_r}},\alpha^{(2)}_{2,n_{k_r}},\alpha^{(3)}_{2,n_{k_r}}, \ldots,\alpha^{(N)}_{2,n_{k_r}},\ldots): r\ge 1\}$ remains dense in $\mathbb{K}^{\mathbb{N}}$.
\end{proof}

B\`es, et al.\!\! \cite{BeMaSa} showed d-hypercyclicity and dense d-hypercyclicity are equivalent for finite families of unilateral weighted backward shifts.  From Corollary \ref{cor:densemanifold} and Corollary \ref{crl:weightedshifts}, we are now able to state that a finite collection of unilateral weighted backward shifts are d-hypercyclic if and only if they posses a dense d-hypercyclic manifold.

In \cite{BeP1}, B\`es and Peris showed that if one considers a finite collection of distinct powers of weighted shifts, d-hypercyclicity is equivalent to satisfying the Disjoint Hypercyclicty Criterion. Of course, the family of powers of weighted shifts is a sub-family of pseudo-shifts and, therefore, it is natural to ask when do finitely many pseudo-shifts satisfy the Disjoint Hypercylicity Criterion. The next theorem tells us precisely when finitely many unilateral pseudo-shifts satisfy the Disjoint Hypercyclicity Criterion.

\begin{theorem} \label{thm:dweakly}
	Let $1 \le p < \infty$ and $N \ge 2$.
    For each integer $i$ with $1 \le i \le N$, let $T_i = T_{f_i, \omega^{(i)}}$ be the unilateral pseudo-shift on $\ell^{p}(\mathbb{N})$ induced by the map $f_i: \mathbb{N} \longrightarrow \mathbb{N}$ and the weight sequence $\omega^{(i)} = (w_{f_i(m)}^{(i)})_{m \in \mathbb{N}}$ as defined in Definition \ref{def:pseudoshift}, and
    let
    $W_{m,n}^{(i)} = \prod_{\nu=1}^{n} w_{f_i^{\nu}(m)}^{(i)}$ as defined in (\ref{intro.3}).
    The following assertions are equivalent:
	\begin{enumerate}
		\item[(i)] The pseudo-shifts $T_1, \ldots, T_N$ are $d$-weakly mixing.
		\item[(ii)] The pseudo-shifts $T_1, \ldots, T_N$ satisfy the Disjoint Hypercyclicity Criterion.
		\item[(iii)] For each integer $R \in \mathbb{N}$, the direct sum operators  $\oplus_{r=1}^{R}T_1,\ldots, \oplus_{r=1}^{R}T_N$ satisfy the Disjoint Blow-Up/Collapse Criterion.
		\item[(iv)] There exists a strictly increasing sequence $(n_k)$ of positive integers which satisfy the following:
		\begin{enumerate}
			\item[(a)] for integers $i,m \in \mathbb{N}$ with $1 \le i \le N$, we have $\left| W_{m,n_k}^{(i)} \right| \longrightarrow  \infty$ as $k \longrightarrow \infty$, and
			\item[(b)] for each $\epsilon>0$ and integers $K, M \in \mathbb{N}$, there exists an integer $k\ge K$ 
			such that for any two integers $i,\ell$ with $1 \le i,\ell \le N$ and $i \neq \ell$, we have $f_{\ell}^{n_k}([M]) \cap f_i^{n_k}([M]) = \emptyset$ and
			\begin{align}
			    \left| 
			    \frac{W_{f_i^{-n_k}(j), n_k}^{(i)}}{W_{f_{\ell}^{-n_k}(j), n_k}^{(\ell)}}
			    \right| < \epsilon
			    \mbox{ \ whenever $j \in f_{\ell}^{n_k}([M]) \cap f_i^{n_k}(\mathbb{N} \setminus [M])$.}
			    \nonumber
			\end{align}
		\end{enumerate}
	\end{enumerate}
\end{theorem}

\begin{proof}
    By Proposition \ref{prop:CoHCrit}, we get (iii) and (ii) are equivalent statements, and (ii) implies (i). It remains to establish (iv) implies (iii) and to establish (i) implies (iv). 
    First, focusing on the establishment of (iv) implies (iii), assume there exists a strictly increasing sequence $(n_k)$ of positive integers which satisfies condition (a) and (b) in statement (iv) of Theorem \ref{thm:dweakly}.  As with the proof of Theorem \ref{thm:characterisation}, we verify the direct sum operators $\oplus_{r=1}^{R}T_1,\ldots, \oplus_{r=1}^{R}T_N$ with $R \in \mathbb{N}$  satisfy the Disjoint Blow-up/Collapse Criterion by utilizing Remark \ref{AlternativeWeakcriterion}.  To this end, let $X_0 = \oplus_{r=1}^{R} \mbox{span} \{ e_m : m \in \mathbb{N} \}$, which is dense in $\oplus_{r=1}^{R} \ell^p(\mathbb{N})$. Consider an arbitrary $\epsilon > 0$, integers $K,L \in \mathbb{N}$ and vectors $y_1, \dots, y_L$ $x_1, \dots, x_N \in X_0$.  
    Write each vector $y_{\ell}$ as $y_{\ell} = (y_{\ell,1}, \dots, y_{\ell,R})$ and each vector $x_i$ as $x_i = (x_{i,1}, \dots, x_{i,R})$ with
	\begin{align}
	    x_{i,r} = \sum_{m=1}^{M} a_{m,r}^{(i)} e_m
	    \mbox{ where 
	    $M = \max \left\{ \mbox{deg}(y_{\ell,r}), \mbox{deg}(x_{i,r}) : 1 \le \ell \le L, 1 \le i \le N, 1 \le r \le R \right\}$.}
	    \label{dweakly.1}
	\end{align}
	For similar reasons as in the proof of Theorem \ref{thm:characterisation}, we may assume the scalars $a_{1,r}^{(i)}, \dots, a_{M,r}^{(i)} \ne 0$ for each integer $i$ with $1 \le i \le N$.  By assumption, we may select an integer $n = n_k$ with $k \ge \max \{ K,M \}$ such that
	\begin{align}
	    \left| W_{m,n}^{(i)} \right|
		& > \frac{M N R \, \Gamma}{\epsilon}
		\mbox{ \ for integers $1 \le i \le N$ and $1 \le m \le M$}
		\label{dweakly.2} \\
		\left| 
		\frac{W_{f_i^{-n}(j), n}^{(i)}}{W_{f_{\ell}^{-n}(j), n}^{(\ell)}} \right|
		& < \frac{\epsilon}{M N R \, \Gamma}
		\mbox{ whenever $i \ne \ell$ and $j \in f_{\ell}^n([M]) \cap f_i^n(\mathbb{N}\setminus [M])$},
		\label{dweakly.3} \\
		 f_{\ell}^n([M]) & \cap f_i^n([M]) = \emptyset 
		\mbox{ whenever $i \ne \ell$},
		\label{dweakly.4}
	\end{align}
	where $\Gamma = \max \{ \Gamma_1, \dots, \Gamma _R \}$ with $\Gamma_r =\max \{ \| x_{i,r} \| : 1 \le i \le N \}$. 
	
	To verify the conditions in Remark \ref{AlternativeWeakcriterion} are satisfied, observe that 
	$T_i^n y_{\ell,r} = 0$ due to the definition of $M$ in (\ref{dweakly.1}) and fact that $n = n_k \ge k \ge \max \{ K,M \}$.  Thus, for integers $i,\ell$ with $1 \le i \le N$ and $1 \le \ell \le L$, we have
    $\left\| ( \oplus_{r=1}^{R} T_i )^n y_{\ell}  \right\|
	= \| (T_i^n y_{\ell,1}, \dots, T_i^n y_{\ell,R}) \| 
	= 0 < \epsilon$.
	Next, for each integer $r$ with $1 \le r \le R$, applying Lemma \ref{lem:technical} with $n=n_k$ and vectors $x_{1,r}, \dots, x_{N,r} \in \mbox{span} \{ e_m : m \in \mathbb{N} \}$ yields a vector $z_r \in \mbox{span} \{ e_m : m \in \mathbb{N} \}$ which satisfies inequalities (\ref{lem:technical.a}), (\ref{lem:technical.b}) and (\ref{lem:technical.c}) in Lemma \ref{lem:technical}.  Consider the vector $z = (z_1, \dots, z_R)$.  For each integer $r$ with $1 \le r \le R$, Lemma \ref{lem:technical} and definition of $\Gamma$ combined with inequality (\ref{dweakly.2}) imply that
	\begin{align}
	    \| z_r \| 
	    & \le M N \, \Gamma_r 
	    \max \left\{ 
	    \left| W_{m,n}^{(\ell)}
	    \right|^{-1} : 1 \le \ell \le N
	    \right\}
	    < M N \, \Gamma_r \frac{\epsilon}{M N R \, \Gamma}
	    \le \frac{\epsilon}{R},
	    \nonumber
	\end{align}
	and so $\| z \| = \| (z_1, \dots, z_R ) \| < \epsilon$.  Lastly, since $f_{\ell}^{n}([M]) \cap f_i^n([M]) = \emptyset$ whenever $i \ne \ell$, it follows from Lemma \ref{lem:technical} and (\ref{dweakly.3}) 
	\begin{align}
	    \| T_i^n z_r - x_{i,r} \| 
	    & \le M N \Gamma_r 
	    \max \left\{ 
	    \left|
	    \frac{W_{f_i^{-n}(j), n}^{(i)}}{W_{f_{\ell}^{-n}(j), n}^{(\ell)}}
	    \right|
	    : \ell \ne i, j \in f_{\ell}^{n}([M]) \cap f_i^n(\mathbb{N} \setminus [M])
	    \right\}
	    < M N \Gamma_r \frac{\epsilon}{M N R \Gamma}
	    \le \frac{\epsilon}{R},
	    \nonumber
	\end{align}
	so 
	$\| (\oplus_{r=1}^{R} T_i)^n z - x_i \| = 
	\| (T_i^n z_1 - x_{i,1}, \dots, T_i^n z_R - x_{i,R} ) \| < \epsilon$ for each integer $i$ with $1 \le i \le N$, which concludes the proof of (iv) implies (iii).
	
	The proof of Theorem \ref{thm:dweakly} is complete after we establish (i) implies (iv).  To this end, select vectors $x = \sum_{m=1}^{\infty} \alpha_m e_m$ and $y = \sum_{m=1}^{\infty} \beta_m e_m$ in $\ell^p(\mathbb{N})$ such that $(x,y)$ is a d-hypercyclic vector for the direct sum operators $T_1 \oplus T_1, \dots, T_N \oplus T_N$.  As in the proof of Theorem \ref{thm:characterisation}, we may select a sequence $(\rho_k)$ of positive scalars for which 
	\begin{align}
	    \rho_k \longrightarrow 0 \mbox{ \ \ and \ \ }
	    \frac{\rho_k}{ \sup \{ | \alpha_m | : m \ge k \} }
	    \longrightarrow \infty
	    \mbox{ \ as $k \longrightarrow \infty$,}
	    \nonumber
	\end{align}
	and select a strictly increasing sequence $(n_k)$ of positive integers such that
	\begin{align}
	    \{ (T_1^{n_k}x, T_1^{n_k}y, \dots, T_N^{n_k}x, T_N^{n_k}y) : k \ge 1 \} 
	    \mbox{ is dense in $\oplus_{i=1}^{2N} \ell^p(\mathbb{N})$, and}
	    \label{dweakly.9} \\
	    \min \left\{ | \langle T_i^{n_k} x, e_m \rangle |: 1 \le i \le N, 1 \le j \le k \right\} > \rho_k 
	    \mbox{ for each integer $k \ge 1$.}
	    \nonumber
	\end{align}
	For the same reason as in the proof of Theorem \ref{thm:characterisation}, it follows that 
	$\left| W_{m,n_k}^{(i)}\right|$ as $k \longrightarrow \infty$.  Next, consider an arbitrary $\epsilon > 0$ and integers $K,M\in \mathbb{N}$. Select a positive constant $\eta$ with  $0 < \frac{10^N \eta}{1-\eta} < \min \{ \epsilon, \frac{1}{10} \}$.  By the density of the set in (\ref{dweakly.9}), there exists an integer $n = n_k$ with $k \ge K$ for which
	\begin{align}
	    \left\| T_i^n x - \sum_{m=1}^{M} e_m \right\| < \eta
	    \mbox{ \ and \ }
	    \left\| T_i^n y - \sum_{m=1}^{M} 10^i e_m \right\| < \eta
	    \mbox{ \ for integers $1 \le i \le N$.}
	    \label{dweakly.11}
	\end{align}
	Similar to the reasoning in the proof of Theorem \ref{thm:characterisation} to generate inequalities (\ref{thm:charact.33}) and (\ref{thm:charact.34}), it follows that
	\begin{align}
	    \left| \alpha_j W_{f_i^{-n}(j)}^{(i)} - 1 \right| 
	    & < \eta, 
	    \mbox{ and }
	    \left| \beta_j W_{f_i^{-n}(j)}^{(i)} - 10^i \right| 
	    < \eta
	    \mbox{ whenever $1 \le i \le N$ and $j \in f_i^n([M])$,}
	    \label{dweakly.12}
	\end{align}
	which in turn implies $\alpha_j, \beta_j \ne 0$ by selection of $\eta$.  Likewise by inequalities in (\ref{dweakly.11}), we also have
	\begin{align}
	    \left| \alpha_j W_{f_i^{-n}(j)}^{(i)} \right| 
	    < \eta
	    \mbox{ whenever $1 \le i \le N$ and $j \in f_i^n(\mathbb{N} \setminus [M])$.}
	    \label{dweakly.13.1}
	\end{align}
	Thus, for integers $i,j,\ell \in \mathbb{N}$ with $i \ne \ell$ and $j \in f_{\ell}^{n}([M]) \cap f_i^{n}(\mathbb{N} \setminus [M])$, we get
	\begin{align}
	   \left| \frac{W_{f_i^{-n}(j)}^{(i)}}{W_{f_{\ell}^{-n}(j)}^{(\ell)}} \right|
	   =  \frac{\left|\alpha_j W_{f_i^{-n}(j)}^{(i)} \right|  }{ \left|\alpha_j W_{f_{\ell}^{-n}(j)}^{(\ell)} \right|  } 
	   < \frac{\eta}{ 1 - \eta} < \epsilon, 
	   \mbox{ by (\ref{dweakly.12}) and (\ref{dweakly.13.1}).}
	   \label{dweakly.14} 
	\end{align}
	Now, if we have  
	$j \in f_{\ell}^{n}([M]) \cap f_i^{n}([M])$, note that
	\begin{align}
	    \left| \frac{W_{f_i^{-n}(j)}^{(i)}}{W_{f_{\ell}^{-n}(j)}^{(\ell)}} - 1 \right|
	    & \le \frac{1}{ \left| \alpha_j W_{f_{\ell}^{-n}(j)}^{(\ell)}  \right| }
	    \left| \alpha_j W_{f_i^{-n}(j)}^{(i)}  - 
	    \alpha_j W_{f_{\ell}^{-n}(j)}^{(\ell)}  \right|
	    \label{dweakly.15} \\
	    & \le \frac{1}{ \left| \alpha_j  W_{f_{\ell}^{-n}(j)}^{(\ell)}  \right| } \left(
	    \left| \alpha_j  W_{f_i^{-n}(j)}^{(i)}  - 1 \right| + \left| 1 - \alpha_j  W_{f_{\ell}^{-n}(j)}^{(\ell)}  \right| 
	    \right)
	    \nonumber \\
	    & < \frac{2 \eta}{1 - \eta}, \mbox{ by (\ref{dweakly.12})}
	    \nonumber \\
	    & < \frac{1}{10}, \mbox{ by selection of $\eta$.}
	    \nonumber
	\end{align}
	By a similar computation, we also have
	\begin{align}
	    \left| \frac{W_{f_i^{-n}(j)}^{(i)}}{W_{f_{\ell}^{-n}(j)}^{(\ell)}} - \frac{10^i}{10^{\ell}} \right|
	    & \le \frac{1}{ \left| \beta_j W_{f_{\ell}^{-n}(j)}^{(\ell)}  \right| }
	    \left| \beta_j W_{f_i^{-n}(j)}^{(i)} - \frac{10^i}{10^{\ell}} \beta_j W_{f_{\ell}^{-n}(j)}^{(\ell)} \right|
	    \label{dweakly.19} \\
	    & \le \frac{1}{ \left| \beta_j W_{f_{\ell}^{-n}(j)}^{(\ell)}  \right| } \left(
	    \left| \beta_j W_{f_i^{-n}(j)}^{(i)}  - 10^i \right| + \frac{10^i}{10^{\ell}}\left| 10^{\ell} - \beta_j W_{f_{\ell}^{-n}(j)}^{(\ell)} \right| \right)
	    \nonumber \\
	    & < \frac{(1 + \frac{10^i}{10^{\ell}}) \eta}{10^{\ell} - \eta}, \mbox{ by (\ref{dweakly.12})}
	    \nonumber \\
	    & < \frac{10^N \eta}{1 - \eta}
	    \nonumber \\
	    & < \frac{1}{10}, \mbox{ by selection of $\eta$.}
	    \nonumber
	\end{align}
	However, both inequalities cannot hold when $i \ne \ell$, and so the set $f_{\ell}^{n}([M]) \cap f_i^{n}([M])$ must be empty whenever $i \ne \ell$.
\end{proof}

By definition, operators $T_1, \dots, T_N$ are d-weakly mixing when the direct sums $T_1 \oplus T_1, \dots, T_N \oplus T_N$ are d-topologically transitive, which is equivalent to the direct sums being densely d-hypercyclic.  The proof of (i) implies (iv) in Theorem \ref{thm:dweakly} only required the direct sums $T_1 \oplus T_1, \dots, T_N \oplus T_N$ of unilateral pseudo-shifts posses a single d-hypercyclic vector $(x,y)$.  From this observation, we conclude that unilateral pseudo-shifts $T_1, \dots, T_N$ are d-weakly mixing if and only if $T_1 \oplus T_1, \dots, T_N \oplus T_N$ are simply d-hypercyclic.  Shkarin, et al.\!\! \cite[Corollary 3.9]{SaSh} showed that in general a direct sums $T_1 \oplus T_1, \dots, T_N \oplus T_N$ of operators being d-hypercycic does not necessarily imply that the operators $T_1, \dots, T_N$ are d-weakly mixing.

A necessary condition for unilateral pseudo-shifts to satisfy the Disjoint Hypercyclicity Criterion is that the underlying maps $f_1, \dots, f_N$ are pairwise distinct.  If there is $f_i \equiv f_{\ell}$ for some $i \ne \ell$, then the condition $f_{\ell}^n([M]) \cap f_i^{n}([M]) = \emptyset$ is never satisfied. To illustrate, B\'es, et al.\!\! \cite{BeMaSa} proved there  unilateral weighted backward shifts $T_1, \dots, T_N$ may be d-hypercyclic, but they never satisfy the Disjoint Hypercyclic Criterion.  Based on Theorem \ref{thm:dweakly}, we can draw the same conclusion by viewing the unilateral weighted backward shifts as unilateral pseudo-shifts as the underlying maps $f_1, \dots, f_n$ are all given by $f_i(m) = m+1$ for all $m \in \mathbb{N}$.

\begin{corollary}\label{cor:unilateralshiftsfail}
Let $T_1, \dots, T_N$ with $N \ge 2$ be unilateral pseudo-shifts on $\ell^{p}(\mathbb{N})$ as defined in Theorem \ref{thm:characterisation}.  The pseudo-shifts $T_1, \dots, T_N$ fail to satisfy the the Disjoint Hypercylicity Criterion if at least two of the underlying maps $f_i, f_{\ell}$ with $i \ne \ell$ are the same map.  In particular, unilateral weighted backward shifts $T_1, \dots, T_N$ always fail to satisfy the Disjoint Hypercyclicity Criterion.
\end{corollary}

 Similar to Corollary \ref{crl:weightedshifts} or Corollary \ref{cor:unilateralshiftsfail}, several existing results for d-hypercyclic shifts can be established using Theorem \ref{thm:characterisation} or Theorem \ref{thm:dweakly}.  For example, B\`es, et al.\!\! \cite[Theorem 4.4]{BeMaSa}  characterized hereditarily d-hypercyclic powers of unilateral weighted shifts; see \cite{BeP1} for definition of hereditarily d-hypercyclicity.

\begin{corollary}\cite[Theorem 4.4]{BeMaSa}\label{cor:dhereditarily}
Let $1 \le p < \infty$ and let $N\ge 2$.  For each integer $i$ with $1 \le i \le N$, let $B_i$ be an unilateral weighted backward shift on $\ell_p(\N)$ with weight sequence $\omega^{(i)} = (w_m^{(i)})$. For distinct integers $1\le r_1<r_2<\cdots<r_N$ and strictly increasing sequence $(n_k)$ of positive integers, the following assertions are equivalent:
	\begin{enumerate}
		\item[(i)] The shifts $B_1^{r_1},B_2^{r_2},\ldots,B_N^{r_N}$ are hereditarily densely d-hypercyclic with respect to the sequence $(n_k)$.
		\item[(ii)] For each $\epsilon>0$ and  integer $M \in \mathbb{N}$, there exists an integer $K \in \mathbb{N}$ such that for each integer $\ell$ with $1 \le \ell \le M$ and each integer $k>K$ the following are satisfied:
		\begin{enumerate}
			\item for each integer $i$ with $1 \le i \le N$, we have \[\left|\prod_{\nu=1}^{r_i n_k}w_{\nu+\ell}^{(i)}\right|>\frac{1}{\varepsilon}\]
			\item for integers $s,t$ with $1 \le s < t \le N$, we have
			\[\left|\prod_{\nu=1}^{r_t n_k}w_{\nu+\ell}^{(t)}\right|>\frac{1}{\varepsilon}\left|\prod_{\nu=1}^{r_s n_k-1}w_{\ell-\nu+r_t n_k}^{(s)}\right|.\]
		\end{enumerate}
	\item[(iii)] The shifts $B_1^{r_1},B_2^{r_2},\ldots,B_N^{r_N}$ satisfy Disjoint Hypercyclicity Criterion with respect to $(n_k)$.
	\end{enumerate}
In particular the shifts $B_1^{r_1},B_2^{r_2},\ldots,B_N^{r_N}$ are d-mixing on $\ell_p(\N)$ if and only if they satisfy the Disjoint Hypercyclicity Criterion with respect to the full sequence $(k)$. 
\end{corollary}

To see how to apply Theorem \ref{thm:dweakly} to establish the above corollary, note that for each integer $i$ with $1 \le i \le N$, we can write $B_i^{r_i} = T_{f_i, \tilde{\omega}^{(i)}}$, where $T_{f_i, \tilde{\omega}^{(i)}}$ is the unilateral lateral pseudo-shift on $\ell^p(\mathbb{N})$ with the map $f_i: \mathbb{N} \longrightarrow \mathbb{N}$ and weight sequence $\tilde{\omega}^{(i)}$ given by 
$f_i(m) = m + r_i$ and $\tilde{w}_{f_i(m)}^{(i)} = w_{m+1}^{(i)}w_{m+2}^{(i)} \cdots w_{m+r_i}^{(i)}$.

\section{Simultaneous Hypercyclic Pseudo-Shifts}\label{SimultaneousPseudoShifts}

While simultaneous hypercyclicity, developed by Bernal and Jung \cite{BGJ}, is a weaker than disjoint hypercyclicity,  simultansous hypercyclicity and disjoint hypercyclic have many parallel results.  In Section \ref{SimultaneousPseudoShifts}, we show a similar situation holds for unilateral pseudo-shits.  We start with the simultaneous version of the blow-up/collapse criterion.

\begin{theorem}[Simultaneous Blow-Up/Collapse Criterion] \label{s-weakcriterion}
	For operators $T_1,\ldots, T_N \in B(X)$ with $N \ge 2$, suppose there exist an increasing sequence $(n_k)$ of positive integers, a dense set $X_0$ in $X$ and maps $S_k:X_0 \to X$ which satisfy the following:
	\begin{enumerate}
		\item[(i)] for each vector $x \in X_0$ and integer $i$ with $1 \le i \le N$, we have $T_i^{n_k}x\longrightarrow 0$ as $k \longrightarrow \infty$, 
		\item[(ii)] for each $\epsilon>0$, integer $K \in \mathbb{N}$ and vector $x_0 \in X_0$, there exists an integer $k\ge K$ such that 
		\begin{enumerate}
			\item[(a)] $\|S_k(x_0)\|<\epsilon$, and
			\item[(b)] $\|T_i^{n_k}S_k(x_0)-x_0\|<\epsilon$ for integers $1 \leq i \leq N$.
		\end{enumerate}
	\end{enumerate}
	Then the operators $T_1, \ldots, T_N$ satisfy the Simultaneous Blow-up/Collapse Property, and hence are densely s-hypercyclic.
\end{theorem}

\begin{proof}
	Let $U,V,W$ be non-empty open sets in $X$ with $0 \in W$.  By the density of the set $X_0$, we may select vectors $x, x_0$ with $x \in V \cap X_0$ and $x_0 \in U \cap X_0$.  By assumptions (i) and (ii) and fact $0 \in W$, there exists an arbitrary large integer $k$ such that $ T_i^{n_k} x \in W$, $S_k(x_0) \in W$ and $T_i^{n_k} S_k(x_0) \in U$ for integers $1 \le i \le N$. 
	This implies $x \in V \cap T_1^{-n_k}(W)\cap \cdots \cap T_N^{-n_k}(W)$ and $S_k(x_0) \in W \cap T_1^{-n_k}(U)\cap \cdots \cap T_N^{-n_k}(U)$, and so
	by Proposition \ref{prop:d-s-BlowUp}, the operators $T_1,\ldots,T_N$ are densely s-hypercyclic.
\end{proof}

The Simultaneous Blow-up/Collapse Criterion has an equivalent formulation.

\begin{remark}\label{rmk:s-alternative}
    The operators $T_1, \dots, T_N \in B(X)$ with $N \ge 2$ satisfy the Simultaneous Blow-up/Collapse Criterion if and only if there is a dense set $X_0$ such that for each $\epsilon > 0$, integers $K,L \in \mathbb{N}$ and vectors $x_0, y_1, \dots, y_L \in X_0$, there exists an integer $n \ge K$ and a vector $z \in X$ with $\| z \| < \epsilon$ such that
    \begin{align}
        \| T_i^{n} y_{\ell} \| < \epsilon 
        \mbox{ \ and \ }
        \| T_i^n z - x_0 \| < \epsilon
        \mbox{ \ for integers $1 \le i \le N$ and $1 \le \ell \le L$.}
        \label{s-alternative.a}
    \end{align}
\end{remark}

Both Theorem \ref{s-weakcriterion} and Remark \ref{rmk:s-alternative} are similar in nature to other sufficient conditions for s-hypercyclicity developed in \cite[Theorem 3.4]{BGJ}.  Using the Simultaneous Blow-up/Collapse Criterion, we establish a complete characterization of s-hypercyclic unilateral pseudo-shifts in terms of their weight sequences.

\begin{theorem} \label{thm:s-characterization}
	Let $1 \le p < \infty$ and $N \ge 2$.
    For each integer $i$ with $1 \le i \le N$, let $T_i = T_{f_i, \omega^{(i)}}$ be the unilateral pseudo-shift on $\ell^{p}(\mathbb{N})$ induced by the map $f_i: \mathbb{N} \longrightarrow \mathbb{N}$ and the weight sequence $\omega^{(i)} = (w_{f_i(m)}^{(i)})_{m \in \mathbb{N}}$ as defined in Definition \ref{def:pseudoshift}, and
    let
    $W_{m,n}^{(i)} = \prod_{\nu=1}^{n} w_{f_i^{\nu}(m)}^{(i)}$ as defined in (\ref{intro.3}).
    The following assertions are equivalent:
    \begin{enumerate}
        \item[(i)] The pseudo-shifts $T_1, \dots, T_N$ are s-hypercyclic.
        \item[(ii)] The pseudo-shifts $T_1, \dots, T_N$ satisfy the Simultaneous Blow-up/Collapse Criterion.
        \item[(iii)] There exists a strictly increasing sequence $(n_k)$ of positive integers which satisfies the following:
        \begin{enumerate}
            \item[(a)] for integers $i,m \in \mathbb{N}$ with $1 \le i \le N$, we have
            $\left| W_{m,n_k}^{(i)} \right| \ \longrightarrow \infty$ as $k \longrightarrow \infty$, and
            \item[(b)] for each $\epsilon > 0$, integers $K, M \in \mathbb{N}$, and  finite collection $\{ a_{m} : 1 \le m \le M \}$ of non-zero scalars in $\mathbb{K} \setminus \{ 0 \}$, there exists an integer $k \ge K$ such that for any integers $i,\ell$ with $1 \le i,\ell \le N$ and $i \ne \ell$, we have
            \begin{align}
                & \left| \frac{W_{f_{i}^{-n_k}(j),n_k}^{(i)}}{W_{f_{\ell}^{-n_k}(j),n_k}^{(\ell)}} \right| < \epsilon 
                \mbox{ \  whenever $j \in f_{\ell}^{n_k}([M]) \cap f_{i}^{n_k}(\mathbb{N} \setminus [M])$,}
                \nonumber
            \end{align}
            and
            \begin{align}
                & \left| \frac{W_{f_{i}^{-n_k}(j),n_k}^{(i)}}{W_{f_{\ell}^{-n_k}(j),n_k}^{(\ell)}} - \frac{a_{f_{i}^{-n_k}(j)}}{a_{f_{\ell}^{-n_k}(j)}}
                \right| < \epsilon 
                \mbox{ \ whenever $j \in f_{\ell}^{n_k}([M]) \cap f_{i}^{n_k}([M])$}.
                \nonumber
            \end{align}
        \end{enumerate}
    \end{enumerate}
\end{theorem}

\begin{proof}
	As with the proofs of the results in Section \ref{DisjointPseudoShfits}, we prove (iii) implies (ii) by assuming there exists a strictly increasing sequence $(n_k)$ of positive integers which satisfy the conditions (a) and (b) in statement (iii), and then use Remark \ref{rmk:s-alternative} with the dense set $X_0 = \mbox{span} \{ e_m : m \in \mathbb{N} \}$ to prove statement (ii).  Consider an arbitrary $\epsilon > 0$, integers $K,L \in \mathbb{N}$ and vectors $x_0, y_1 \dots, y_L \in X_0$.  Write the vector $x_0$ as
	\begin{align}
	    x_0 = \sum_{m=1}^{M} a_m e_m 
	    \mbox{ \ where  
	    $M = \max \{ \mbox{deg}(x_0), \mbox{deg}(y_1), \dots, \mbox{deg}(y_L) \}$.}
	    \label{s-charact.1}
	\end{align}
	We may further assume $a_1, \dots, a_M \ne 0$.  
	By assumption, there exists an integer $n = n_k$ with $k \ge \max \{ K, M \}$ such that
	\begin{align}
	    \left| W_{m,n}^{(i)} \right| 
	    & > \frac{M N \Gamma}{\epsilon}
	    \mbox{ for integers $1 \le i \le N$ and $1 \le m \le M$,}
	    \label{s-charact.2} \\
	    \left| \frac{W_{f_i^{-n}(j),n}^{(i)}}{W_{f_{\ell}^{-n}(j),n}^{(\ell)}} \right|
	    & < \frac{\epsilon}{2 M N\Gamma }
	    \mbox{ whenever $i\ne \ell$ and $j \in f_{\ell}^n([M]) \cap f_i^n(\mathbb{N} \setminus [M])$, and}
	    \label{s-charact.3} \\
	    \left| \frac{W_{f_i^{-n}(j),n}^{(i)}}{W_{f_{\ell}^{-n}(j),n}^{(\ell)}} - \frac{a_{f_i^{-n}(j)}}{a_{f_{\ell}^{-n}(j)}} \right|
	    & < \frac{\epsilon}{2 M N \Gamma }
	    \mbox{ whenever $i\ne \ell$ and 
	    $j \in f_{\ell}^n([M]) \cap f_i^n([M])$,}
	    \label{s-charact.4}
	\end{align}
	where $\Gamma = \| x_0 \|$.  Lemma \ref{lem:technical} applied with integers $n, M$ and vectors $\overbrace{x_0, \dots, x_0}^{N-times}$ provides a vector $z \in X_0$ with
	\begin{align}
	    \| z \| 
	    & \le M N \Gamma 
	    \max \left\{ \left| W_{m,n}^{(i)} \right|^{-1} : 1 \le i \le N, 1 \le m \le M
	    \right\}
	    \label{s-charact.5} \\
	    \| T_1^n z - x_0 \|
	    & \le M N \Gamma \max \left\{
	    \left| \frac{W_{f_1^{-n}(j),n}^{(1)}}{W_{f_{\ell}^{-n}(j),n}^{(\ell)}} \right| : 2 \le \ell \le N, j \in f_{\ell}^n([M]) \cap f_1^n(\mathbb{N} \setminus [M])
	    \right\},
	    \label{s-charact.6}
	   \end{align}
and for each integer $i$ with $2 \le i \le N$,
	   \begin{align}
	    \| T_i^n z - x_0 \|
	    & \le M N \Gamma 
	    \max \left\{ \left| \frac{W_{f_i^{-n}(j),n}^{(i)}}{W_{f_{\ell}^{-n}(j),n}^{(\ell)}}
	    - \frac{a_{f_i^{-n}(j)}}{a_{f_{\ell}^{-n}(j)}} \right| : \ell \ne i, j \in f_{\ell}^n([M]) \cap f_i^n([M])
	    \right\}
	    \label{s-charact.7} \\
	    & \ \ \ \ +
	     M N \Gamma 
	    \max \left\{ \left| \frac{W_{f_i^{-n}(j),n}^{(i)}}{W_{f_{\ell}^{-n}(j),n}^{(\ell)}}
	     \right| : \ell \ne i, j \in f_{\ell}^n([M]) \cap f_i^n(\mathbb{N} \setminus [M])
	    \right\}.
	    \nonumber
	\end{align}
	It follows from inequalities (\ref{s-charact.2}), (\ref{s-charact.3}), (\ref{s-charact.4}) with (\ref{s-charact.5}), (\ref{s-charact.6}), (\ref{s-charact.7}) that $\| z \| < \epsilon$ and $\|T_i^n z - x_0 \| < \epsilon$ for integers $1 \le i \le N$ as required in Remark \ref{rmk:s-alternative}.  Furthermore, by the definition of the integer $M$ in (\ref{s-charact.1}) and the selection of $n = n_k \ge \max \{ K, M \}$,
	 we get $\| T_i^n y_{\ell} \| = 0 < \epsilon$ for integers $i, \ell$ with $1 \le i \le N$ and $1 \le \ell \le L$, which completes the establishment of the conditions in Remark \ref{rmk:s-alternative}.
	
	Since Theorem \ref{s-weakcriterion} establishes (ii) implies (i), we complete the proof of Theorem \ref{thm:s-characterization} by proving (i) implies (iii). Let 
	$x = \sum_{m=1}^{\infty} \alpha_m e_m$ be an s-hypercyclic vector for the unilateral pseudo-shifts $T_1, \dots, T_N$.  As in the proof of Theorem \ref{thm:characterisation}, since $x$ is a s-hypercyclic vector, there exists a sequence $(\rho_k)$ of positive scalars such that $\rho_k \longrightarrow 0$ and $\frac{\rho_k}{\sup \{ |\alpha_m| : m \ge k \} } \longrightarrow \infty$ as $k \longrightarrow \infty$ and a strictly increasing sequence $(n_k)$ of positive integers such that 
	\begin{align}
	    & \{ (T_1^{n_k} x, \dots, T_N^{n_k} x) : k \ge 1 \}
	    \mbox{ is dense in $\Delta(\oplus_{i=1}^{N}\ell^p(\mathbb{N}))$ and}
	    \nonumber \\
	    & \min \left\{ 
	    \left| \langle T_i^{n_k}x, e_j \rangle \right| :
	    1 \le i \le N, 1 \le j \le k \right\} > \rho_k
	    \mbox{\ for integers $k \ge 1$}.
	    \label{s-charact.9}
	\end{align}
	Following the same reasoning as in (\ref{thm:charact.31}) in the proof of Theorem \ref{thm:characterisation}, we get 
	$\left| W_{m,n_k}^{(i)} \right| \longrightarrow \infty$ as $k \longrightarrow \infty$, which establishes condition (a) in statement (iii).  
	
	To establish condition (b), consider any arbitrary $\epsilon > 0$, integers $K,M \in \mathbb{N}$ and finite collection $\{ a_m : 1 \le m \le M \}$ of non-zero scalars in $\mathbb{K}\setminus \{ 0 \}$.  
	Select positive constants $C, \eta$  such that
	\begin{align}
	    1 < C \min \left\{ | a_m | : 1 \le m \le M 
	    \right\}, 
	    \ \ 
	    \max \left\{  \left|
	    \frac{a_j}{a_m} \right| : 1 \le j,m \le M \right\}
	    < C
	    \mbox{ \ and \ }
	    0 < \frac{(1+C) \eta}{1- \eta} < \epsilon.
	    \nonumber
	\end{align}
	By the density of the set in (\ref{s-charact.9}), there exists an integer $k \ge K$ such that
	\begin{align}
	    \left\| T_i^{n_k} x - \sum_{m=1}^{M} C e_m \right\| < \eta
	    \mbox{ \ for integers $1 \le i \le N$.}
	    \nonumber
	\end{align}
	Following the same arguments as in (\ref{thm:charact.33}), (\ref{thm:charact.34}), (\ref{thm:charact.35}) and (\ref{thm:charact.36}) in the proof of Theorem \ref{thm:characterisation}, we get
	\begin{align}
	    \left| \frac{W_{f_{i}^{-n_k}(j),n_k}^{(i)}}{W_{f_{\ell}^{-n_k}(j),n_k}^{(\ell)}} \right| 
	    &  < \epsilon 
        \mbox{ \  whenever $j \in f_{\ell}^{n_k}([M]) \cap f_{i}^{n_k}(\mathbb{N} \setminus [M])$, and }
        \nonumber \\
        \left| \frac{W_{f_{i}^{-n_k}(j),n_k}^{(i)}}{W_{f_{\ell}^{-n_k}(j),n_k}^{(\ell)}} - \frac{a_{f_{i}^{-n_k}(j)}}{a_{f_{\ell}^{-n_k}(j)}} \right| 
        & < \epsilon 
        \mbox{ \  whenever $j \in f_{\ell}^{n_k}([M]) \cap f_{i}^{n_k}([M])$,}
        \nonumber
	\end{align}
	which concludes the proof of (i) implies (iii).
\end{proof}

The Simultaneous Blow-up/Collapse Criterion provides a sufficient condition for operators $T_1, \dots, T_N$ to be densely s-hypercyclic.  Thus, from Theorem \ref{thm:s-characterization}, we get s-hypercycicity and dense s-hypercyclicity are equivalent ideas with unilateral pseudo-shifts.

\begin{corollary}\label{cor:s-dense}
Let $T_1, \dots, T_N$ with $N \ge 2$ be unilateral pseudo-shifts on $\ell^p(\mathbb{N})$ as defined in Theorem \ref{thm:s-characterization}.  The pseudo-shifts $T_1, \dots, T_N$ are s-hypercyclic if and only if they are densely s-hypercyclic.
\end{corollary}

The Simultaneous Hypercyclicity Criterion is an alternative sufficient condition for establishing dense s-hypercyclicity for a finite family of operators.  Similar to the disjoint hypercyclicity situation, when focusing on unilateral pseudo-shifts, satisfying the Simultaneous Hypercyclicity Criterion has an equivalent characterization expressed in terms of the weight sequences. 

\begin{theorem}\label{thm:s-HyperCriterionCharacterization}
	Let $1 \le p < \infty$ and $N \ge 2$.
    For each integer $i$ with $1 \le i \le N$, let $T_i = T_{f_i, \omega^{(i)}}$ be the unilateral pseudo-shift on $\ell^{p}(\mathbb{N})$ induced by the map $f_i: \mathbb{N} \longrightarrow \mathbb{N}$ and the weight sequence $\omega^{(i)} = (w_{f_i(m)}^{(i)})_{m \in \mathbb{N}}$ as defined in Definition \ref{def:pseudoshift}, and
    let
    $W_{m,n}^{(i)} = \prod_{\nu=1}^{n} w_{f_i^{\nu}(m)}^{(i)}$ as defined in (\ref{intro.3}).
    The following assertions are equivalent:
	\begin{enumerate}
		\item[(i)] The pseudo-shifts $T_1,\ldots,T_N$ are s-weakly mixing.
		\item[(ii)] The pseudo-shifts $T_1,\ldots,T_N$ satisfy the Simultaneous Hypercyclicity Criterion.
		\item[(iii)] For each integer $R \in \mathbb{N}$, the direct sum operators $\oplus_{r=1}^R T_1,\ldots,\oplus_{r=1}^R T_N$ satisfy the Simultaneous Blow-Up/Collapse Criterion.
		\item[(iv)] There exists a strictly increasing sequence $(n_k)$ of positive integers which satisfy the following: 
		\begin{enumerate}
			\item[(a)] for integers $i, m \in \mathbb{N}$ with $1 \le i \le N$, we have $\left| W_{m,n_k}^{(i)} \right| \longrightarrow \infty$ as $k \rightarrow \infty$, and
			\item[(b)] for each $\epsilon>0$ and integers $K,M \in \mathbb{N}$, there exists an integer $k \ge K$ such that for any two integers $i, \ell$ with  $1 \le i, \ell \le N$ and  $i \neq \ell$, we have
			\begin{align}
				\left| \frac{W_{f_i^{-n_k}(j),n_k}^{(i)}}{W_{f_{\ell}^{-n_k}(j),n_k}^{(\ell)}}
				\right| < \epsilon
				\mbox{ \ whenever $j \in f_{\ell}^{n_k}([M]) \cap f_i^{n_k}(\mathbb{N} \setminus [M])$,}
				\nonumber 
			\end{align}
		    and for any integer $j \in f_{\ell}^{n_k}([M]) \cap f_i^{n_k}([M])$, we have
		    $f_{\ell}^{-n_k}(j) = f_i^{-n_k}(j)$ and 
		    \begin{align}
		        \left| \frac{W_{f_i^{-n_k}(j),n_k}^{(i)}}{W_{f_{\ell}^{-n_k}(j),n_k}^{(\ell)}} - 1
		        \right| < \epsilon.
		        \nonumber
		    \end{align}
		\end{enumerate}
	\end{enumerate}
\end{theorem}

\begin{proof}
	The proof Theorem \ref{thm:s-HyperCriterionCharacterization} is similar to the proof of Theorem \ref{thm:dweakly}.  From Proposition \ref{prop:sHC}, we get (iii) and (ii) are equivalent statements, and (ii) implies (i). Thus, it remains to show (iv) implies (iii) and show (i) implies (iv). 
	
	Starting with the proof of (iv) implies (iii), suppose the the strictly increasing sequence $(n_k)$ satisfies conditions (a) and (b) in statement (iv).  With Remark \ref{rmk:s-alternative} in mind, let $X_0 = \oplus_{r=1}^{R} \mbox{span} \{ e_m : m \in \mathbb{N} \}$ with $R \in \mathbb{N}$, and
	consider an arbitrary $\epsilon > 0$, integers $K,L \in \mathbb{N}$ and vectors $x_0, y_1 \dots, y_L \in X_0$.  Write each vector $y_{\ell}$ as $y_{\ell} = (y_{\ell,1}, \dots, y_{\ell,R})$ and write $x_0$ as $x_0 = (x_{0,1}, \dots, x_{0,R})$ with
	\begin{align}
	    x_{0,r} = \sum_{m=1}^{M} a_{m,r} e_m 
	    \mbox{ where 
	    $M = \max \{ \mbox{deg}(x_{0,r}), \mbox{deg}(y_{\ell,r}) : 
	    1 \le i \le N, 1 \le \ell \le L, 1 \le r \le R \}$.}
	    \label{s-Hyper.1}
	\end{align}
	Further assume $a_{1,r}, \dots, a_{M,r} \ne 0$ for integer each integer $r$ with $1 \le r \le R$.  
	By assumption, there exists an integer $n = n_k$ with $k \ge \max \{ K, M \}$ such that
	\begin{align}
	    \left| W_{m,n}^{(i)} \right| 
	    & > \frac{R M N \Gamma}{\epsilon}
	    \mbox{ for integers $1 \le i \le N$ and $1 \le m \le M$,}
	    \label{s-Hyper.2} \\
	    \left| \frac{W_{f_i^{-n}(j),n}^{(i)}}{W_{f_{\ell}^{-n}(j),n}^{(\ell)}} \right|
	    & < \frac{\epsilon}{2 R M N \Gamma }
	    \mbox{ whenever $i \ne \ell$ and $j \in f_{\ell}^n([M]) \cap f_i^n(\mathbb{N} \setminus [M])$, and}
	    \label{s-Hyper.3} \\
	    \left| \frac{W_{f_i^{-n}(j),n}^{(i)}}{W_{f_{\ell}^{-n}(j),n}^{(\ell)}} - 1 \right|
	    & < \frac{\epsilon}{2 R M N \Gamma }
	    \mbox{ and $f_{\ell}^{-n}(j) = f_i^{-n}(j)$ whenever $i\ne \ell$ and $j \in f_{\ell}^n([M]) \cap f_i^n([M])$}
	    \label{s-Hyper.4}
	\end{align}
	where $\Gamma = \max \{ \Gamma_1, \dots, \Gamma_r \}$ with $\Gamma_r = \| x_{0,r} \|$.  Following the same arguments in the proofs of the previous results, the selection of $n,M$ implies that 
	$\| (\oplus_{r=1}^{R}T_i)^n y_{\ell}  \| = 
	\| (T_i^n y_{\ell,1}, \dots, T_i^n y_{\ell, R} ) \| = 0 < \epsilon$ for integers $i,\ell$ with $1 \le i \le N$ and $1 \le \ell \le L$.  Next,
	or each integer $r$ with $1 \le r \le R$, Lemma \ref{lem:technical} applied with integers $n, M$ and vectors $\overbrace{x_{0,r}, \dots, x_{0,r}}^{N-times}$ yields a vector $z_r$ with
	\begin{align}
	    \| z_r \| 
	    & \le M N \Gamma_r 
	    \max \left\{ \left| W_{m,n}^{(i)} \right|^{-1} : 1 \le i \le N, 1 \le m \le M
	    \right\}
	    \label{s-Hyper.5} \\
	    \| T_1^n z_r - x_{0,r} \|
	    & \le M N \Gamma_r \max \left\{
	    \left| \frac{W_{f_1^{-n}(j),n}^{(1)}}{W_{f_{\ell}^{-n}(j),n}^{(\ell)}} \right| : 2 \le \ell \le N, j \in f_{\ell}^n([M]) \cap f_i^n(\mathbb{N} \setminus [M])
	    \right\},
	    \label{s-Hyper.6} 
	   \end{align}
    and for each integer $i$ with $2 \le i \le N$,
	   \begin{align}
	    \| T_i^n z_r - x_{0,r} \|
	    & \le M N \Gamma_r 
	    \max \left\{ \left| \frac{W_{f_i^{-n}(j),n}^{(i)}}{W_{f_{\ell}^{-n}(j),n}^{(\ell)}}
	    - \frac{a_{f_i^{-n}(j)}}{a_{f_{\ell}^{-n}(j)}} \right| : \ell \ne i, j \in f_{\ell}^n([M]) \cap f_i^n([M])
	    \right\}
	    \label{s-Hyper.7} \\
	    & \ \ \ \ +
	     M N \Gamma_r 
	    \max \left\{ \left| \frac{W_{f_i^{-n}(j),n}^{(i)}}{W_{f_{\ell}^{-n}(j),n}^{(\ell)}}
	     \right| : \ell \ne i, j \in f_{\ell}^n([M]) \cap f_i^n(\mathbb{N} \setminus [M])
	    \right\}.
	    \nonumber
	\end{align}
	Consider the vector $z = (z_1, \dots, z_r)$. Note that inequalities (\ref{s-Hyper.2}), (\ref{s-Hyper.3}) with (\ref{s-Hyper.5}), (\ref{s-Hyper.6}) and the  definition of $\Gamma$ yield $\| z_r \| < \frac{\epsilon}{R}$ and $
	\| T_1^n z_r - x_{0,r} \| < \frac{\epsilon}{R}$, which in turn implies that
	$ \| z \| = \| (z_1, \dots, z_r ) \| < \epsilon$ and 
	$\| (\oplus_{r=1}^{R} T_1)^n z - x_0 \| = \| (T_1^n z_1 - x_{0,1}, \dots, T_1^n z_r - x_{0,r}) \| < \epsilon$.  Furthermore,  whenever $j \in f_{\ell}^{n}([M]) \cap f_i^{n}([M])$ with $i \ne \ell$, we have by assumption that $f_{\ell}^{-n}(j) = f_i^{-n}(j)$, which in turn implies
	\begin{align}
	    \left| \frac{W_{f_i^{-n}(j),n}^{(i)}}{W_{f_{\ell}^{-n}(j),n}^{(\ell)}}
	    - \frac{a_{f_i^{-n}(j)}}{a_{f_{\ell}^{-n}(j)}} \right|
	    =
	    \left| \frac{W_{f_i^{-n}(j),n}^{(i)}}{W_{f_{\ell}^{-n}(j),n}^{(\ell)}}
	    - 1 \right|.
	    \nonumber
	\end{align}
	The above observation together with inequalities (\ref{s-Hyper.3}), (\ref{s-Hyper.4}) and (\ref{s-Hyper.7}) gives us $\| T_i^n z_r - x_{0,r} \| < \frac{\epsilon}{R}$ for each integer $i$ with $2 \le i \le N$, and so 
	$\| (\oplus_{r=1}^{R}T_i)^n z - x_0 \| 
	= \| (T_i^n z_1 - x_{0,1}, \dots, T_i^n z_R - x_{0,R}) \|
	< \epsilon$ for $2 \le i \le N$.  
	
	Finally, to prove (i) implies (iv), select vectors $x = \sum_{m=1}^{\infty} \alpha_m e_m $ and $y = \sum_{m=1}^{\infty} \beta_m e_m$ such that $(x,y)$ is an s-hypercyclic vector for the direct sum operators $T_1 \oplus T_1, \dots, T_N \oplus T_N$.  As in the previous arguments, we may select a sequence $(\rho_k)$ of positive scalars such that $\rho_k \longrightarrow 0$ and 
	$\frac{\rho_k}{\sup \{ |\alpha_m| : m \ge k \} } \longrightarrow \infty$ as $k \longrightarrow \infty$, and select a strictly increase sequence $(n_k)$ of positive integers such that
	\begin{align}
	    & \{ (T_1^{n_k} x, T_1^{n_k} y, \dots, T_N^{n_k} x, T_N^{n_k} y ) : k \ge 1 \}
	    \mbox{ is dense in $\Delta(\oplus_{i=1}^{2N} \ell^p(\mathbb{N}))$, and }
	    \label{s-Hyper.10} \\
	    & \min \left\{ \left| \langle T_i^{n_k}x, e_m \rangle \right| : 1 \le i \le N, 1 \le j \le k \right\} > \rho_k 
	    \mbox{ for integers $k \ge 1$.} 
	    \nonumber
	\end{align}
	As seen in the previous proofs, the above conditions establishes condition (a) in statement (iv).  To establish condition (b), consider an arbitrary $\epsilon > 0$ and integers $K, M \in \mathbb{N}$.  Select a positive constant $\eta$ with $0 < \frac{10^N \eta}{1 - \eta} < \min\{ \epsilon, \frac{1}{10} \}$.  By the density of the set in (\ref{s-Hyper.10}), we may select an integer $n = n_k$ with $k \ge K$ for which
	\begin{align}
	    \left\|  T_i^n x - \sum_{m=1}^{N} e_m \right\| < \eta
	    \mbox{ and }
	    \left\|  T_i^n y - \sum_{m=1}^{N} 10^m e_m \right\| < \eta
	    \mbox{ for integers $1 \le i \le N$.}
	    \label{s-Hyper.12}
	\end{align}
	Similar to the reasoning in the proof of Theorem \ref{thm:characterisation} to generate inequalities (\ref{thm:charact.33}) and (\ref{thm:charact.34}), it follows from (\ref{s-Hyper.12}) that
	\begin{align}
	    \left|\alpha_j W_{f_i^{-n}(j)}^{(i)}  - 1 \right| 
	    < \eta
	    \mbox{ and }
	    \left| \beta_j W_{f_i^{-n}(j)}^{(i)} - 10^{f_i^{-n}(j)} \right| 
	    < \eta
	    \mbox{ whenever $1 \le i \le N$ and $j \in f_i^n([M])$}
	    \nonumber
	   \end{align}
	   and
	   \begin{align}
	    & \left| \alpha_j W_{f_i^{-n}(j)}^{(i)} \right| 
	    < \eta
	    \mbox{ whenever $1 \le i \le N$ and $j \in f_i^n(\mathbb{N} \setminus [M])$.}
	    \nonumber
	\end{align}
	Following similar computations used to establish inequalities (\ref{dweakly.14}), (\ref{dweakly.15}) and (\ref{dweakly.19}) in the proof of Theorem \ref{thm:dweakly}, the above inequalities yield
	\begin{align}
	    \left| \frac{W_{f_i^{-n}(j)}^{(i)}}{W_{f_{\ell}^{-n}(j)}^{(\ell)}} \right| < \frac{\eta}{1 - \eta} < \epsilon
	    \mbox{ whenever $i \ne \ell $ and $j \in f_{\ell}^n([M]) \cap f_i^{n}(\mathbb{N} \setminus [M])$.}
	    \nonumber
	\end{align}
	Also, whenever $j \in f_{\ell}^n([M]) \cap f_i^{n}([M])$ with $i \ne \ell$, we have 
	\begin{align}
	    \left| \frac{W_{f_i^{-n}(j)}^{(i)}}{W_{f_{\ell}^{-n}(j)}^{(\ell)}} - 1 \right| 
	    & < \frac{2 \eta}{1 - \eta} < \min \left\{ \epsilon, \frac{1}{10} \right\} \le \epsilon
	    \mbox{ \ and \ }
	    \left| \frac{W_{f_i^{-n}(j)}^{(i)}}{W_{f_{\ell}^{-n}(j)}^{(\ell)}} - \frac{10^{f_i^{-n}(j)}}{10^{f_{\ell}^{-n}(j)}} \right| 
	    < \frac{10^N \eta}{1 - \eta} < \frac{1}{10}.
	    \nonumber
	\end{align}
	Both above inequalities simultaneously occur only when $f_i^{-n}{(j)} = f_{\ell}^{-n}{(j)}$, which completes the proof of condition (b) in statement (iv).
\end{proof}

If the pseudo-shifts $T_1, \dots, T_N$ with $N \ge 2$ satisfy the Simultaneous Hypercyclicity Criterion, then by Theorem \ref{thm:s-characterization} and Theorem \ref{thm:s-HyperCriterionCharacterization}, it follows that those same pseudo-shifts satisfy the Simultaneous Blow-up/Collapse Criterion.  However, the following example shows the converse fails to hold. That is, the Simultaneous Blow-up/Collapse Criterion is a weaker statement.

\begin{example}\label{ex:s-BlowUpCollapse}
Let $1 \le p < \infty$ and let $N \ge 2$.  For each  constant $\beta > 1$, there exist  unilateral pseudo-shifts $T_1, \dots, T_N$ on $\ell^p(\mathbb{N})$ with $\| T_1 \| = \cdots = \| T_N \| = \beta$ that satisfy the Simultaneous Blow-up/Collapse Criterion and yet fail to satisfy the Simultaneous Hypercyclicity Criterion.
\end{example}

\begin{proof}
    First we construct two of the desired unilateral pseudo-shifts $T_1, T_2$ on $\ell^p(\mathbb{N})$, and then extend the example for $N \ge 3$ unilateral pseudo-shifts $T_1, \dots, T_N$.  Begin by letting $(\widetilde{w}_m^{(1)})_{m \in \mathbb{N}}, (\widetilde{w}_m^{(2)})_{m \in \mathbb{N}}$ be the respective weight sequences for d-hypercyclic unilateral weighted backward shifts $B_1, B_2$ on $\ell^p(\mathbb{N})$ with
    \begin{align}
        \| B_i \| = \sup \left\{ | \widetilde{w}_m^{(i)} | :
        m \in \mathbb{N}
        \right\}
        = \beta \mbox{ for integers $i=1,2$.}
        \nonumber
    \end{align}
    See \cite[Proposition 2.3]{BeMaSa} for the existence of two such shifts $B_1,B_2$.  Consider the unilateral pseudo-shifts $T_1 = T_{f_1, \omega^{(1)}}$ and $T_2 = T_{f_2, \omega^{(2)}}$ associated with the strictly increasing maps $f_1, f_2: \mathbb{N} \longrightarrow \mathbb{N}$ given by
    \begin{align}
        f_1(m) =
        \begin{cases}
            4, & \text{ if $m=1$} \\
            5, & \text{ if $m=2$} \\
            2m, & \text{ otherwise} 
        \end{cases}
        \mbox{ \ and \ }
        f_2(m) =
        \begin{cases}
            2m, & \text{ if $m=2^r$ for some $r \ge 0$} \\
            2m+1, & \text{ otherwise} 
        \end{cases},
        \nonumber
    \end{align}
    and weight sequences $\omega^{(1)}, \omega^{(2)}$ given by
    \begin{align}
        w_{m}^{(i)} =
        \begin{cases}
            \widetilde{w}_r^{(i)}, & \text{ if $m = 2^r$ for some $r \ge 0$} \\
            \beta, & \text{ otherwise}
        \end{cases}.
        \label{s-example.3}
    \end{align}
    Note that $\| T_i \| = \sup \left\{ | w_{f_i(m)}^{(i)} | : m \in \mathbb{N} \right\} = \beta$ for integers $i=1,2$.
    
    To verify the pseudo-shifts $T_1, T_2$ fail to satisfy the Simultaneous Hypercyclicity Criterion, observe that for each integer $\nu \in \mathbb{N}$, we have
    \begin{align}
        f_1^{\nu}(m) =
        \begin{cases}
            2^{\nu + 1}, & \text{ if $m=1$} \\
            5\cdot 2^{\nu - 1}, & \text{ if $m=2$} \\
            m2^{\nu}, & \text{ otherwise} 
        \end{cases}
        \mbox{ \ and \ }
        f_2^{\nu}(m) =
        \begin{cases}
            2^{\nu + r}, & \text{ if $m=2^r$ for some $r \ge 0$} \\
            \text{odd integer}, & \text{ otherwise} 
        \end{cases}.
        \label{s-example.4}
    \end{align}
    In particular, note that $f_1^{\nu}(m)$ is an even integer for any $m, \nu \in \mathbb{N}$ with $\nu \ge 2$ and $f_2^{\nu}(m)$ is only an even integer when $f_2^{\nu}(m)$ is a power of 2.  Thus, the sets $f_{\ell}^n([M]) \cap f_{i}^n([M])$ and 
    $f_{\ell}^n([M]) \cap f_{i}^n(\mathbb{N} \setminus [M])$ with $i \ne \ell$ consist only of powers of 2.  More specifically, for each integer $M \in \mathbb{N}$, select $r_M \in \mathbb{N}$ such that $2^{r_M} \le M < 2^{1+ r_M}$, and so for each integer $n \ge 2$, we have
    \begin{align}
        f_{1}^n([M]) \cap f_{2}^n(\mathbb{N} \setminus [M])
        & = \{ 2^{n+1}, 2^2 2^n, \dots, 2^{r_M} 2^n \}
        \cap \{ 2^{n+1+r_M}, 2^{n+2+r_M}, \dots \} = \emptyset,
        \nonumber \\
        f_{2}^n([M]) \cap f_{1}^n(\mathbb{N} \setminus [M])
        & = \{ 2^n, 2^{n+1}, 2^{n+2}, \dots, 2^{n+r_M} \}
        \cap \{ 2^{1+r_M}2^n, 2^{2+r_M}2^n, \dots \} = \emptyset, 
        \mbox{ and }
        \nonumber \\
        f_{1}^n([M]) \cap f_{2}^n([M])
        & = \{ 2^{n+1}, 2^{n+2}, \dots, 2^{n+r_M} \}.
        \nonumber
    \end{align}
    Observe that $2^{n+1} \in f_{1}^n([M]) \cap f_{2}^n([M])$ for any  $M,n \in \mathbb{N}$ with $n \ge 2$, but
    $f_1^{-n}(2^{n+1}) = 1 \ne 2 = f_2^{-n}(2^{n+1})$ by (\ref{s-example.4}), which in turn implies condition (b) in statement (iv) in Theorem \ref{thm:s-HyperCriterionCharacterization} is never satisfied.  Hence, the unilateral pseudo-shifts $T_1, T_2$ fail to satisfy the Simultaneous Hypercyclicity Criterion.
    
    We verify the pseudo-shifts $T_1, T_2$ satisfy the Simultaneous Blow-up/Collapse Criterion by establishing statement (iii) in Theorem \ref{thm:s-characterization}.  To this end, let 
    $\{ \gamma_k : k \in \mathbb{N} \}$ be a countable dense set in $\mathbb{K} \setminus \{ 0 \}$.  Since the weight sequence $(\widetilde{w}_m^{(1)})_{m \in \mathbb{N}}, (\widetilde{w}_m^{(2)})_{m \in \mathbb{N}}$ are associated with d-hypercyclic unilateral weighted backward shifts, by Theorem \ref{thm:characterisation}, we may select a strictly increasing sequence sequence $(n_k)$ of positive integers such that for each integer $k \ge 1$, we have
    \begin{align}
        & \beta^{n_k} > k \mbox{ \ and \ } 
        \left| \prod_{\nu=1}^{n_k} \widetilde{w}_{\nu + r}^{(i)} \right| > k
        \mbox{ for integers $i=1,2$ and $1 \le r \le k$,}
        \label{s-example.8} \\
        & \left| \prod_{\nu=1}^{n_k} 
        \frac{\widetilde{w}_{\nu + 1}^{(1)}}{\widetilde{w}_{\nu + 1}^{(2)}} - \gamma_k \right| < \frac{1}{k}
        \mbox{ \ and \ }
        \left| \prod_{\nu=1}^{n_k} 
        \frac{\widetilde{w}_{\nu + 1}^{(2)}}{\widetilde{w}_{\nu + 1}^{(1)}} - \frac{1}{\gamma_k} \right| < \frac{1}{k},
        \label{s-example.9} \\
        & \left| \prod_{\nu=1}^{n_k} 
        \frac{\widetilde{w}_{\nu + r}^{(1)}}{\widetilde{w}_{\nu + r}^{(2)}} - 1 \right| < \frac{1}{k}
        \mbox{ \ and \ }
        \left| \prod_{\nu=1}^{n_k} 
        \frac{\widetilde{w}_{\nu + r}^{(2)}}{\widetilde{w}_{\nu + r}^{(1)}} - 1 \right| < \frac{1}{k}
        \mbox{ \ for integer $2 \le r \le k$.}
        \label{s-example.10}
    \end{align}
    To show we have the desired sequence $(n_k)$,  observe that for integers $m,n \in \mathbb{N}$ with $n \ge 2$, by (\ref{s-example.3}) and (\ref{s-example.4}), we have
    \begin{align}
        W_{m,n}^{(1)} = \prod_{\nu=1}^{n_k} w_{f_1^{\nu}(m)}^{(1)} 
        & =
        \begin{cases}
            \prod_{\nu=1}^{n_k} w_{2^{\nu+1}}^{(1)}, 
            & \text{ if $m = 1$} \\
            \prod_{\nu=1}^{n_k} w_{5 \cdot 2^{\nu-1}}^{(1)}, 
            & \text{ if $m = 2$} \\
            \prod_{\nu=1}^{n_k} w_{m2^{\nu}}^{(1)}, 
            & \text{ otherwise}
        \end{cases}
        \label{s-example.11} \\
        & =
        \begin{cases}
            \prod_{\nu=1}^{n_k} \widetilde{w}_{\nu+1}^{(1)}, 
            & \text{ if $m = 1$} \\
            \prod_{\nu=1}^{n_k} \widetilde{w}_{\nu+r}^{(1)}, 
            & \text{ if $m = 2^r$ for some $r \ge 2$} \\
            \beta^{n_k}, 
            & \text{ otherwise}
        \end{cases}.
        \nonumber 
    \end{align}
    For similar reasons, we have
    \begin{align}
         W_{m,n}^{(2)} = \prod_{\nu=1}^{n_k} w_{f_2^{\nu}(m)}^{(2)} =
         \begin{cases}
            \prod_{\nu=1}^{n_k} \widetilde{w}_{\nu+r}^{(2)}, 
            & \text{ if $m = 2^r$ for some $r \ge 0$} \\
            \beta^{n_k}, 
            & \text{ otherwise}
            \label{s-example.12}
        \end{cases}.
    \end{align}
    By (\ref{s-example.8}), we get $\left|  W_{m,n}^{(i)}\right| \longrightarrow \infty$ as $k \longrightarrow \infty$, establishing condition (a) in statement (iii) in Theorem \ref{thm:s-characterization}.  For establishing condition (b) in the same statement, consider an arbitrary $\epsilon > 0$, integers $K, M \in \mathbb{N}$ and finite set
    $\{ a_1, \dots, a_M \}$ of nonzero scalars in $\mathbb{K}\setminus \{ 0 \}$.  
    By the density of the set 
    $\{ \gamma_k : k \in \mathbb{N} \}$ in $\mathbb{K} \setminus \{ 0 \}$, we may select an integer $k > \max \{ K, r_M \}$, where $2^{r_M} \le M < 2^{1+ r_M}$, such that
    \begin{align}
        \left| \frac{a_1}{a_2} - \gamma_k \right| 
        < \frac{\epsilon}{2}, \quad
        \left| \frac{a_2}{a_1} - \frac{1}{\gamma_k} \right| 
        < \frac{\epsilon}{2}
        \mbox{ \ and \ }
        \frac{1}{k} < \frac{\epsilon}{2}.
        \label{s-example.13}
    \end{align}
    Since $f_{\ell}^{n_k}([M]) \cap f_{i}^{n_k}(\mathbb{N} \setminus [M]) = \emptyset$ for $i \ne \ell$, we only examine 
    $j \in f_{1}^{n_k}([M]) \cap f_{2}^{n_k}([M]) = \{ 2^{n_k+1}, \dots, 2^{n_k+r_M} \}$.
    For $j = 2^{n_k+1} \in f_{1}^{n_k}([M]) \cap f_{2}^{n_k}([M])$, we have
    \begin{align}
         \left| 
        \frac{W_{f_1^{-n_k}(2^{n_k+1}), n_k}^{(1)}}{W_{f_2^{-n_k}(2^{n_k+1}), n_k}^{(2)}}
         - \frac{a_{f_1^{-n_k}(2^{n_k+1})}}{a_{f_2^{-n_k}(2^{n_k+1})}}
         \right|
        & = 
        \left| 
        \frac{W_{1, n_k}^{(1)}}{W_{2, n_k}^{(2)}}
        - \frac{a_{1}}{a_{2}}
        \right|, 
        \mbox{ by (\ref{s-example.4})}
        \nonumber \\
        & \le
        \left| 
        \prod_{\nu=1}^{n_k} \frac{\widetilde{w}_{\nu + 1}^{(1)}}{\widetilde{w}_{\nu + 1}^{(2)}}
        - \gamma_k
        \right|
        + 
        \left| \gamma_k - \frac{a_{1}}{a_{2}} \right|,
        \mbox{ by (\ref{s-example.11}), (\ref{s-example.12})}
        \nonumber \\
        &< \frac{1}{k} + \frac{\epsilon}{2} < \epsilon, 
        \mbox{ by (\ref{s-example.9}), (\ref{s-example.13}).}
        \nonumber 
    \end{align}
    Likewise, we have
    \begin{align}
        \left| \frac{W_{f_2^{-n_k}(2^{n_k+1}), n_k}^{(2)}}{W_{f_1^{-n_k}(2^{n_k+1}), n_k}^{(1)}}
        - \frac{a_{f_2^{-n_k}(2^{n_k+1})}}{a_{f_1^{-n_k}(2^{n_k+1})}}
        \right| < \epsilon.  
        \nonumber
    \end{align}
    For $j = 2^{n_k+r} \in f_{1}^{n_k}([M]) \cap f_{2}^{n_k}([M])$ with $2 \le r \le r_M$, we have
    \begin{align}
         \left| 
        \frac{W_{f_1^{-n_k}(2^{n_k+r}), n_k}^{(1)}}{W_{f_2^{-n_k}(2^{n_k+r}), n_k}^{(2)}}
         - \frac{a_{f_1^{-n_k}(2^{n_k+r})}}{a_{f_2^{-n_k}(2^{n_k+r})}}
         \right|
        & = 
        \left| 
        \frac{W_{2^r, n_k}^{(1)}}{W_{2^r, n_k}^{(2)}}
        - \frac{a_{2^r}}{a_{2^r}}
        \right|, 
        \mbox{ by (\ref{s-example.4})}
        \nonumber \\
        & =
        \left| 
        \prod_{\nu=1}^{n_k} \frac{\widetilde{w}_{\nu + r}^{(1)}}{\widetilde{w}_{\nu + r}^{(2)}}
        - 1
        \right|
        \mbox{ by (\ref{s-example.11}), (\ref{s-example.12})}
        \nonumber \\
        &< \frac{1}{k}  < \epsilon, 
        \mbox{ by (\ref{s-example.10}), (\ref{s-example.13}).}
        \nonumber 
    \end{align}
    For similar reasons, we have
    \begin{align}
        \left| 
        \frac{W_{f_2^{-n_k}(2^{n_k+r}), n_k}^{(2)}}{W_{f_1^{-n_k}(2^{n_k+r}), n_k}^{(1)}}
         - \frac{a_{f_2^{-n_k}(2^{n_k+r})}}{a_{f_1^{-n_k}(2^{n_k+r})}}
         \right| < \epsilon.
         \nonumber 
    \end{align}
    Therefore, the pseudo-shifts $T_1, T_2$ satisfy the Simultaneous Blow-up/Collapse Criterion by Theorem \ref{thm:s-characterization}.
    
    We conclude our example by noting that we can extend this example to $N \ge 3$ unilateral pseudo-shifts by taking 
    $T_1, T_2, \overbrace{T_1, \dots, T_1}^{N-2 \ times}$.  This collection of unilateral pseudo-shifts remains s-hypercyclic (and so satisfy the Simultaneous Blow-up/Collapse Criterion), but fail to satisfy the Simultaneous Hypercyclicity Criterion as desired.
\end{proof}

Bernal-Gonz\'alez and Jung in \cite{BGJ} characterized s-hypercyclicity of distinct powers of weighted shifts. Complementing this result, we give a full characterization for s-hypercyclic unilateral weighted shifts. Contrary to the disjoint hypercyclicity situation, we show that finitely many weighted shifts are s-hypercyclic if and only if they satisfy the  Simultaneous Hypercyclicity Criterion.

\begin{corollary} \label{cor:simultaneousweightedshifts}
Let $1 \le p < \infty$ and let $N \ge 2$.  For each integer $i$ with $1 \le i \le N$, let $T_i$ be a unilateral weighted backward shift on $\ell^p(\mathbb{N})$ with the weight sequence $\omega^{(i)} = (w^{(i)}_m)_{m \in \mathbb{N}}$.   
The following assertions are equivalent:

\begin{enumerate}
	\item[(i)] The shifts $T_1, \ldots, T_N$ are $s$-hypercyclic.
	\item[(ii)] The shifts $T_1, \ldots, T_N$ satisfy the Simultaneous Hypercyclicity Criterion.
	\item[(iii)] There exists a strictly increasing $(n_k)$ of positive integers which satisfy the following:
	\begin{enumerate}
	    \item[(a)] for integers $i,m \in \mathbb{N}$ with $1 \le i \le N$, we have 
	    $\left| \prod_{\nu=1}^{n_k} w_{m + \nu}^{(i)} \right| \longrightarrow \infty $ as $k \longrightarrow \infty$, and
	    \item[(b)] for each $\epsilon > 0$ and integers $K,M \in \mathbb{N}$, there exists an integer $k \ge K$ such that
	    \begin{align}
	        \left| \prod_{\nu = 1}^{n_k} 
	        \frac{w_{m+\nu}^{(i)}}{w_{m+\nu}^{(\ell)}} - 1
	        \right| < \epsilon
	        \mbox{ for integers $1 \le i,\ell \le N$ and $1 \le m \le M$.}
	        \nonumber
	    \end{align}
	\end{enumerate}
\end{enumerate}
\end{corollary}

\begin{proof}
    As noted in the proof of Corollary \ref{crl:weightedshifts}, the unilateral weighted shift $T_i$ with weight sequence $\omega^{(i)} = (w_{m}^{(i)})_{m \in \mathbb{N}}$ can be expressed as the unilateral pseudo-shift $T_{f_i, \omega^{(i)}}$ where $f_i(m) = m+1$ for every $m \in \mathbb{N}$.  Furthermore, for integers $i, \ell, m \in \mathbb{N}$ with $1 \le i, \ell \le N$, we have
    \begin{align}
        W_{m,n}^{(i)} = \prod_{\nu = 1}^{n} w_{m+\nu}^{(i)}
        \mbox{ \ and \ }
        \frac{W_{m,n}^{(i)}}{W_{m,n}^{(\ell)}} = 
        \prod_{\nu = 1}^{n} \frac{w_{m+\nu}^{(i)}}{w_{m+\nu}^{(\ell)}}.
        \nonumber
    \end{align}
    The result now follows from Theorem \ref{thm:s-HyperCriterionCharacterization} and the observation that 
    $f_{\ell}^n([M]) \cap f_i^n(\mathbb{N} \setminus [M]) = \emptyset$ whenever $i \ne \ell$,  
    and $f_i^{-n}(j) = j - n = f_{\ell}^{-n}(j)$ for each integer $j   \in \{ n+1, \dots, n+M \} = f_{\ell}^n([M]) \cap f_i^n([M])$.
\end{proof}

Salas \cite{Sa} showed that a direct sum $T_1 \oplus \cdots \oplus T_N$ of unilateral weighted backwards shifts is hypercyclic if and only if their corresponding weight sequences $\omega^{(1)}, \dots, \omega^{(N)}$ satisfy 
\begin{align}
    \sup \left\{
    \min \left\{ \left| \prod_{\nu = 1}^{n} w_{1 + \nu}^{(i)} \right| : 1 \le i \le N \right\}: n \ge 1\right\} = \infty.
    \label{DirectSum.1}
\end{align}
By combining Salas' result with condition (a) in statement (iii) of Corollary \ref{cor:simultaneousweightedshifts}, it immediately follows that the direct sum $T_1 \oplus \cdots \oplus T_N$ is hypercyclic whenever $T_1, \dots, T_N$ are s-hypercyclic unilateral weighted backward shifts.

\begin{corollary}\label{cor:DirectSumShifts}
    Let $T_1, \dots, T_N$ with $N \ge 2$ be unilateral weighted backward shifts on $\ell^p(\mathbb{N})$ as defined in Corollary \ref{cor:simultaneousweightedshifts}.  If the shifts $T_1, \dots, T_N$ are s-hypercyclic, then the direct sum 
    $T_1 \oplus \cdots \oplus T_N$ is hypercyclic.
\end{corollary}

We conclude our section with a simple example showing the converse of Corollary \ref{cor:DirectSumShifts} fails hold.  Select two distinct constants $\alpha, \beta > 1$, and let $T_1, T_2$ be the unilateral weighted backard shifts with the respective weight sequences $\omega^{(1)} $, $\omega^{(2)}$ given by $w_m^{(1)} = \alpha$ for each $m \in \mathbb{N}$,  $w_2^{(2)} = \beta$ and $w_m^{(2)} = \alpha$ otherwise.  Since $\alpha, \beta > 1$, it follows from (\ref{DirectSum.1}) that the direct sum $T_1 \oplus T_2$ is hypercyclic.  However, the shifts $T_1, T_2$ fail to be s-hypercyclic.  To verify, note that for any integer $n \in \mathbb{N}$, we have
\begin{align}
    \left| \frac{W_{1,n}^{(1)}}{W_{1,n}^{(2)}} - 1 \right| 
    = \left| \frac{\alpha^n}{\beta \alpha^{n-1}} - 1 \right|
    = \left| \frac{\alpha}{\beta} - 1 \right|,
    \nonumber 
\end{align}
and so condition (b) in statement (iii) of Corollary \ref{cor:simultaneousweightedshifts} is not satisfied for any $\epsilon$ with $0 < \epsilon < | \frac{\alpha}{\beta} -1 |$.





\end{document}